\documentclass[times,final]{elsarticle}
\usepackage[utf8]{inputenc}
\usepackage[english]{babel}
\usepackage{amssymb}
\usepackage{latexsym}
\usepackage{amsmath}
\usepackage{mathrsfs}
\usepackage{url}
\usepackage{xcolor}

\title{Fast method to simulate dynamics of two-phase medium with intense interaction between phases by smoothed particle hydrodynamics: gas-dust mixture with polydisperse particles, linear drag, one-dimensional tests}

\author{Olga Stoyanovskaya$^{1,2}$,  Maxim Davydov$^{1,2}$, Tamara Markelova$^{1,2,3}$, Maxim Arendarenko$^{1,2}$, Elizaveta Isaenko$^2$,  Valeriy Snytnikov$^{1,2,3}$ }

\address{$^1$Lavrentiev Institute of Hydrodynamics SB RAS, Ave. Lavrentieva, 15, Novosibirsk, 630090, Russia}
\address{$^2$Novosibirsk State University, Str. Pirogova, 1, Novosibirsk, 630090, Russia}
\address{$^3$Boreskov Institute of Catalysis SB RAS, Ave. Lavrentieva, 5, Novosibirsk, 630090, Russia}
\ead{o.p.sklyar@gmail.com}

\begin{document}

\maketitle

\begin{abstract}
To simulate the dynamics of fluid with polydisperse particles on macroscale level, one has to solve hydrodynamic equations with several relaxation terms, representing momentum transfer from fluid to particles and vice versa. For small particles, velocity relaxation time (stopping time) can be much shorter than  dynamical time of fluid \textcolor{red}{that makes} this problem stiff and thus computationally expensive. 
We present a new fast method for computing several stiff drag terms in two-phase polydisperse medium with Smoothed Particle Hydrodynamics (SPH). In our implementation, fluid and every fraction of dispersed phase are simulated with different sets of particles. The method is based on (1) linear interpolation of velocity values in drag terms, (2) implicit approximation of drag terms that conserves momentum with machine precision\textcolor{red}{,} and (3) solution of system of $N$ linear algebraic equations with $O(N^2)$ arithmetic operation instead of $O(N^3)$.

We studied the properties of the proposed method on one-dimensional problems with known solutions. We found that we can obtain acceptable accuracy of the results with numerical resolution independent of short stopping time values. All simulation results discussed in the paper are obtained with \textcolor{red}{open source} software. 

\end{abstract}

\section{Introduction}

%%%%%%%%%%%%%%%%%%%%%%%%%%

\subsection{Definition and examples of two-phase polydisperse media}

Two-phase polydisperse medium is a mixture of carrier fluid (compressible or incompressible) and dispersed inclusions of other phases (solid grains, bubbles or drops), in which features of dispersed phase differ by several orders of magnitude. For example, dusty gas with particles of different sizes is a two-phase polydisperse medium. In general, we can talk about the distribution of elements of the dispersed phase with respect to their features: size, mass, velocity, etc.

Two-phase polydisperse media are found in many technological processes and devices. These are combustion chambers of heat engines, air ducts of solid propellant and liquid jet engines, fire extinguishing systems, devices for sand and shot blasting of various surfaces, pneumatic conveyors of loose materials \cite{pneumo2010}, dust collectors (aspirators) of various types, chemical reactors with fluid catalyst \cite{fluidization2013}, technique of cold and detonation spraying of particles on surface, etc \cite{Varaksin2013}. Aerodisperse flows are studied in the design of aerospace craft, which move at supersonic speed in the atmosphere of the Earth, Mars and other planets \cite{Vasilevskii2001}.

Two-phase polydisperse flows also occur in a large number of natural phenomena. Tornadoes, sandstorms, ashfalls in the atmosphere caused by volcanic eruptions, forest fires, rains, clouds and many other processes are well known. It is of great interest to model the dynamics of gas and dust medium in astrophysics studying the problems of origin of planets in circumstellar disks \cite{HaworthEtAl2016} and chemical evolution in molecular clouds. 

\subsection{Approaches to description of two-phase polydisperse flows}

The extensive literature is devoted to mathematical models of laminar and turbulent flows of a fluid with dispersed inclusions of another phase, for example, \cite{Balachandar2010,Soo1989,Nigmatullin,DeichFilippovBook,Gidaspow1994,Babukha1967,Marble1970,Sternin}. A two-phase medium can be described at three levels (see for example, \cite{marchisio_fox_2013}). 

The first of these levels is a microscale level with direct modeling of the dynamics of fluid molecules and particles of dispersed inclusions. This approach requires the greatest computational costs. 

The second level is a macroscale level, where the dynamics of fluid and dispersed inclusions are described by hydrodynamic approximation. On this level, we derive the governing equations from original microscale equations via averaging with respect to volume or ensemble of particles.

This approach requires significantly less computational resources compared to the costs in molecular dynamics. In the hydrodynamic approach, the carrier phase is described by a system of hydrodynamic equations, and the dispersed phase is represented by a set of $N$ monodisperse fractions. For each of the fractions, the continuity equations, equations of motion and energy are written implying that each fraction moves at its own velocity (and has its own temperature). The result is a mathematical model of multifluid hydrodynamics. The examples of such models in engineering and astrophysics can be found in \cite{Ivanenko2019,Tukmakova2019,Lozhkin2014,Korolev2018,Drajkowska2019}. The main feature of this approach is that in every local volume the velocity  \textbf{of each fraction} of the dispersed phase has a single value, and there is no particle trajectory crossing (PTC) \textbf{inside of one fraction}. Another example of macroscopic model is one-fluid approach, in which the equations for the total density and barycentric velocity of gas-dust medium and the diffusion equation for a dispersed phase are solved. For a gas-polydisperse dust medium, the derivation of one-fluid equations from multi-fluid hydrodynamics can be found in \cite{LaibePrice2014OneFluidDust}.

The third is a meso-scale approach, in which the carrier is described as a fluid, and the dynamics of dispersed particles is described by using kinetic equations. The third approach is more costly than the second one, and contains the second approach as an extreme special case \cite{marchisio_fox_2013,UGKS2019}. The third approach admits different particle velocity values in one local volume.

In the paper, we use a macroscopic multifluid approach to the description of the flows of a fluid with polydisperse particles (see governing equations (\ref{eq:contSPHgas})-(\ref{eq:motionSPHdust}) in section~\ref{sec:equations}). This approach is valid if the mean free path of molecules for the carrier medium and particles of each fraction of the dispersed phase is much smaller than the length scale of the system. Generally, this condition is expressed in terms of the dimensionless Knudsen number. The mean free path is inversely proportional to the number density of elements in the medium. Thus, the condition can be violated if the dispersed phase has a low concentration. However, there is no need to refuse the approximation of a fluid for low concentrations of the dispersed phase if the stopping time of particles (velocity relaxation time or time during which a particle velocity reaches a stationary value relative to a gas velocity) is much shorter than the dynamical time of the problem. This condition is expressed in terms of the Stokes number (see details in \cite{marchisio_fox_2013}). 

\subsection{Momentum exchange between phases - stiff relaxation terms, computational problems}

In fluid with polydisperse inclusions, carrier and dispersed phases can exchange mass, momentum and energy. Momentum and energy exchange is described by velocity relaxation terms in equations of motion and thermal relaxation terms in equations of energy. The velocity relaxation term is proportional to the relative velocity between phases and inversely proportional to stopping time (see e.g. (\ref{eq:motionSPHgas})-(\ref{eq:motionSPHdust}) terms with $t_i$ in denominator). The presence of these terms can make the problem stiff if relaxation time values are much less than dynamical time of the carrier phase. In this case, stopping time is a small parameter of the problem.

Problems with stiff relaxation terms are actively studied in applied \cite{StiffRelaxationTerm} and computational \cite{Jin2010} mathematics. In addition to dynamics of multiphase media, such problems have many applications (plasma physics \cite{PlasmaAP2017}, elasticity with memory, etc.). The numerical solution of problems with stiff relaxation terms is computationally challenging. 

It is necessary to use a time step less than stopping time for application of explicit schemes to such tasks. This constraint on time step for two-phase polydisperse media can be much stiffer than the Courant condition for pure fluids. To avoid such a time step constraint due to stability condition, the stiff relaxation terms are approximated implicitly. However, many authors showed \cite{NumericsStiff} that even in this case, in order to obtain acceptable accuracy of the solution, there is still need to adjust spatial and time step according to the value of a small parameter. In particular, the authors in \cite{StoyanovskayaDust} considered the example of particle moving under the gas drag and gravitation of a massive body. They showed how implicit approximation of drag together with operator splitting with respect to physical processes can significantly reduce the accuracy of the solution. In this example, the authors considered approximation when gas affects the particle velocity through the drag force, and particles do not affect the gas velocity. Such a problem contained one stiff relaxation term. The computing of drag becomes more complex when the influence of all $N$ fractions of the disperse phase on carrier is considered, then $2N$ relaxation terms are in the system. For $N=1$, the examples of schemes with implicit approximation of drag and with different accuracy can be found in \cite{StoyanovskayaDustMultigrain,StoyanovskayaNonlinear,Vorobyov2018}.

In some cases, you can omit finding a numerical solution to the original problem and replace it with a simplified problem derived from the original one via expansion with respect to a small parameter (for example, \cite{StiffRelaxationTerm}). In particular, for the dynamics of gas-dust medium, this means that the velocity of the dispersed phase is found not from the solution of the original differential equation, but from the algebraic relation \cite{Marble1970}. A simplified differential-algebraic problem becomes nonstiff and simpler from a computational point of view. However, the transition to it is possible only for the extreme case when stopping time is much less than dynamical time of carrier phase and does not allow modeling of transient cases with stopping time close to dynamical time.

%%%%%%%%%%%%%%%%%%%%%%%%%%%%%%%%%%%%%%

\subsection{Fluids with mono- and polydisperse inclusions - approximation of stiff relaxation terms in grid methods}

The practice of application of implicit approximations of stiff relaxation terms showed that in many cases, if we keep the numerical resolution fixed, the numerical scheme approximation error increases as small parameters decrease. This is due to the fact that the main term of the error is proportional to $(\displaystyle\frac{\tau}{t_{\rm stop}})^{\alpha}$, where $\tau$ is a time step, $t_{\rm stop}$ is stopping time (details can be found in Table 3 of work \cite{StoyanovskayaDust} and formulas (17)-(18) in \cite{StoyanovskayaDustMultigrain}). However, in some cases, the opposite is observed: as a small parameter decreases, the error of the numerical solution also decreases. This is due to the fact that the main term of the error is multiplied by the value, which at $t_{\rm stop} \longrightarrow 0$ tends to zero faster than $(\displaystyle \frac{\tau}{t_{\rm stop}})^{\alpha}$ increases. Such schemes have been called asymptotic preserving (AP) schemes. 

To obtain a scheme that preserves the asymptotic (providing a bounded solution with acceptable accuracy for $\tau \gg t_{\rm stop}$), we have to approximate stiff relaxation terms using ideas described in \cite{Jin2010,StoyanovskayaDustMultigrain}. The first idea is to approximate stiff relaxation terms implicitly, i.e., to use values from the next time moment in the computations of the terms. The second idea is to solve analytically the Cauchy problem for a linear system of ordinary differential equations at each step and to get the solution at $t+\tau$. For the dynamics of fluid with mono- and polydisperse inclusions, the approaches to approximation of stiff relaxation terms in grid methods are presented in \cite{Miniati2010,StoyanovskayaDustMultigrain,Li2017,BenitezLlambay2019}. All of them are based on the idea of implicit approximation of drag ensuring the implementation of the law of momentum conservation locally. That is, AP methods ensure that in each computational cell the momentum lost by particles due to gas drag is equal to the momentum acquired by gas due to dust drag. The violation of the local momentum conservation law in grid methods leads to a bias of the wave propagation velocity in the medium (see illustrations in \cite{StoyanovskayaDustMultigrain,Vorobyov2018}). The transfer of momentum between the phases depends on the relative velocity between the carrier and the dispersed phases. Therefore, it is easy to provide the implementation of the local momentum conservation law in the Euler methods if the computational grids for a fluid and dispersed phases coincide. However, it is non-trivial to transfer this idea to Lagrangian methods.  

\subsection{Fluids with mono- and polydispersed inclusions - approximation of stiff relaxation terms in Smoothed Particle Hydrodynamics}

One of common Lagrangian methods for modeling of fluid dynamics is the Smoothed Particle Hydrodynamics. In this method, a fluid is replaced by a set of particles that are carriers of mass, momentum and energy of the medium. These particles are the moving nodes of approximation. A smooth approximation function is constructed to calculate spatial derivatives using values known in nodes. The differentiation operation is applied to the smooth approximant. In multi-fluid smoothed particle hydrodynamics, it is assumed that the carrier phase and each fraction of the dispersed phase are modeled by their own set of particles. Obviously, the computational costs and the asymptotic properties of the method will depend on the way the relative velocities in the relaxation terms are calculated. To date, three ways of spatial velocity interpolation have been developed for the smoothed particle hydrodynamics, which we will call particle-particle \cite{MonaghanKocharyan1995}, duplicate node \cite{BateDust2014} and drag in a cell \cite{IDIC2018}.
 
The particle-particle method is proposed by Monaghan and Kocharyan \cite{MonaghanKocharyan1995}. In this method, neighboring particles of the dispersed phase act on each particle of the carrier phase (and vice versa), in this case, the contribution of each neighbor depends on the distance between the particles and is determined by the kernel of the method. 

The implicit approximation of the drag force by the method of Monaghan and Kocharyan \cite{MonaghanKocharyan1995} is nontrivial because each particle requires using the radius-vectors and velocity values of all its neighbors at the next time moment. However, several algorithms \cite{Monaghan1997,Monaghan2020,LaibePrice2011Test,LaibePriceAstroDrag} have been developed so far to reduce the computational cost of each time step.

The Monaghan and Kocharyan method guarantees the conservation of momentum and angular momentum in the whole computational domain. However, for an arbitrary volume inside the domain, the fulfillment of the momentum conservation law is not guaranteed. As a result, in the calculation of intense drag, a numerical overdissipation appears. This overdissipation was demonstrated and discussed in papers \cite{LaibePrice2011Test,BateDust2014,IDIC2018}. Laibe and Price \cite{LaibePrice2011Test} showed that the scale of this overdissipation depends on stopping time. To avoid this, we have to use smoothing length constrained by
\begin{equation}
\label{eq:h}
    h<t_{\rm stop} c_{\rm s},
\end{equation}
where $c_{\rm s}$ is a sound speed in pure gas \cite{LaibePrice2011Test}. The authors \cite{IDIC2018} showed that the longer the kernel is, the greater the overdissipation scale is, and moreover, the scale may exceed $h$ by an order of magnitude.

Due to numerical overdissipation problem in the Monaghan-Kocharyan method another way of drag term computing - a duplicate node - was proposed and tested in \cite{BateDust2014,BateDust2015,Booth2015}. In this way, we compute the velocity of dust at the point where the gas particle is located (and vice versa) using standard SPH interpolation. The duplicate node method preserves the Lagrangian nature of the interphase interaction computing and reduces the numerical overdissipation of the solution in the case of stiff drag \cite{2018CompareMethods,2019CompareMethods}, but does not guarantee global and local conservation of momentum and angular momentum. 

This motivated the authors \cite{IDIC2018} to propose the Lagrangian-Eulerian approach to compute stiff mutual  interaction between phases (drag in a cell), while retaining the computing of other forces completely Lagrangian. To do this, they proposed to split all the particles into non-crossing cells and calculate the average arithmetic velocity of phases in each cell. Then, one has to compute the velocity of individual particles of the carrier phase using the cell-averaged velocity values of the dispersed phase (and vice versa), so that the law of momentum conservation is satisfied in each cell. This approach of velocity interpolation is similar to implementation of NGP kernel in Particle-in-Cell method, which recently has been applied to astrophysical simulations of dusty gas with stiff mutual drag \cite{YangJohansen2016,YangJohCarrera2017}. Studying this method, they found \cite{IDIC2018} that in the computation of intense interphase interaction, this approach allows us to use spatial resolution (smoothing length) and cell size independent of a small parameter, i.e., stopping time.

Thus, all the numerical algorithms of two-fluid smoothed particle hydrodynamics presented in literature can be classified according to the method of approximation of stiff relaxation term in time and the method of spatial velocity interpolation. The resulting classification is shown in Table \ref{tab:resume}.

A separate branch of the Lagrangian methods for computing the dynamics of multiphase media is the SPH approximation of the equations of one-fluid hydrodynamics \cite{LaibePriceOneFluidSPH}. In this case, each particle carries the parameters of the carrier and dispersed phases and the need for velocity interpolation disappears, which significantly reduces the computational cost. The one-fluid approach guarantees local (and hence global) conservation of momentum in the medium, but does not guarantee the conservation of mass of the dispersed phase. 

By now, the one-fluid SPH approach has already been extended to the case of a polydisperse medium \cite{Multigrain}. In this paper, we firstly extend the multi-fluid SPH approach to polydisperse media using "implicit drag in a cell" way to compute stiff relaxation term. Owing to the Lagrangian nature of the SPH method, when approximating equations of gas-dynamic type, the condition of uniqueness of velocity in the local volume of the medium can be violated and the trajectories of model particles can cross. For a guaranteed "continuity" of the medium, artificial viscosity must be introduced into the calculation. Therefore, the idea of computing drag for multifluid hydrodynamics proposed in the paper does not base on the "single-velocity" assumption of the medium. It means that the idea can be promising for the implementation of mathematical models of meso-level where kinetic equations for the disperse phase are solved by particle-based method. For the mesoscale approach, AP schemes have already been developed for monodisperse mixtures \cite{CARRILLO2008, GOUDON2013145}.

\subsection{Structure of the paper}
\label{sec:intr_structure}
Section~\ref{sec:equations} provides the governing equations of gas-polydisperse dust dynamics and the formulation of two one-dimensional test problems: the propagation of acoustic waves and a shock tube in polydisperse dusty gas. Section~\ref{sec:NumMethod} shows the main idea of new way of intense drag computing and formulas of the numerical algorithm for the whole problem with detailed description of "drag-in-cell" formulas in subsection~\ref{sec:IDICMultigrain}. To study the properties of the new algorithm, we implemented it in the form of code in C-language, which we publish at

 \verb|https://github.com/MultiGrainSPH/1D_Dust_DS.git| as a supplementary material for the paper. 
Section~\ref{sec:results} provides the results of simulation of test problems in the case of small and multi-scale stopping time values. The conclusions are given in Section~\ref{sec:resume}. In \ref{sec:MK}, we describe our implementation of the Monaghan-Kocharyan method that we used in comparative studies. \ref{sec:type} shows that the isothermal gas-polydisperse dust system being solved is hyperbolic. \ref{sec:dustyWaveMultigrain} describes the way of obtaining a reference solution for the first test problem, i.e., sound waves in an isothermal gas-polydisperse dust medium, and also describes the Scilab code we published at

\verb|https://github.com/MultiGrainSPH/1D_Dust_DS/tree/master/DustyWave|, which generates this reference solution. \ref{sec:analysis} describes the way of obtaining a reference solution for the shock-tube, the second test problem, namely, the effective sound speed is derived for a polydisperse gas-dust medium in the extreme case of small stopping time values of all dispersed phases.

\begin{table*}
\centering
\begin{minipage}{160mm}
\caption{Methods for computing stiff drag in a two-phase monodisperse medium for two-fluid smoothed particle hydrodynamics (TFSPH). In the work marked by $^*$, a particular case is considered where the dynamics of the dispersed phase does not affect the dynamics of the carrier. See also recent work \cite{PriceLaibe2020} on TFSPH published after the first submission of our paper.}
  \begin{center}
   \label{tab:resume}
  \begin{tabular}{cccc}
  \hline
    Time \ Space & Particle-particle (kernel) & Duplicate node & Cell drag\\
   \hline
   Implicit      &  Monaghan 1997 \cite{Monaghan1997} &   &   \\
   approximation &  Monaghan 2020 \cite{Monaghan2020} &   &  Stoyanovskaya \\
                 &  Laibe, Price 2012 \cite{LaibePrice2011Test,LaibePriceAstroDrag} &   & et al. 2018 \cite{IDIC2018} \\
 \hline
  "Exact" solution        &   & \textcolor{magenta}{Loren-Aguilar}   &  \\
        &   & Bate 2014,2015 \cite{BateDust2014,BateDust2015} &   \\     
   & & Booth et al. 2015$^*$ \cite{Booth2015}& \\
  
\end{tabular}
\end{center}
\end{minipage}
\end{table*}  

\section{Governing equations}

\label{sec:equations}

We consider a two-phase polydisperse medium, in which the carrier phase is a compressible non-viscous gas of density $\rho$, velocity $\mathbf{v}$ and internal energy $e$. The gas contains $N$ fractions of dispersed inclusions differing in particle size. For each fraction, the condition of applicability of hydrodynamic medium description is satisfied and the following space-averaged values are defined: $\mathbf{u}_i$ is velocity, $\rho_i$ is mass density, $t_{i}$ is stopping time. The particles of all fractions consist of the same substance of material density $\rho_{\rm s}$, have a spherical shape, and occupy a zero volume in space. The gas-dust medium has total pressure, which does not act on the particles due to a zero-volume fraction of particles. Mass and heat exchange between gas and particles is absent. All dust fractions exchange momentum with gas but not with each other. The momentum exchange occurs in the mode of Epstein (Knudsen) or Stokes drag. This means that one of two conditions holds. The first of them is 
\begin{equation}
\label{eq:Knudsen}
   \max_{i=1}^N {s_i} < 2.25 \lambda, 
\end{equation} where $s_i$ is the radius of fraction particle $i$, $\lambda$ is the length of free path of the carrier gas molecules. The second equation is
\begin{equation}
\label{eq:Stokes}
    \max_{i=1}^N Re_{i \rm drag}<1,
\end{equation}
where $Re_{i \rm drag}=4\displaystyle \frac{\| \mathbf{v} -\mathbf{u}_i \|}{c_{\rm s}} \frac{s_i}{\lambda}$ is the Reynolds number calculated from the relative velocity between gas and particles. When satisfying (\ref{eq:Knudsen}) or (\ref{eq:Stokes}), $t_{i}$ is independent of $\mathbf{v}$ and $\mathbf{u}_i$. 

We write the equations of continuity and motion for gas and dust fractions in terms of a total spatial derivative

\begin{equation}
\label{eq:contSPHgas}
    \displaystyle\frac{ \partial \rho_{\rm g}}{\partial t}+ ( \mathbf{v} \cdot \nabla \rho_{\rm g} ) = -\rho_{\rm g} \nabla \mathbf{v},
\end{equation}

\begin{equation}
\label{eq:contSPHdust}
    \displaystyle\frac{ \partial \rho_i}{\partial t}+( \mathbf{u}_i \cdot \nabla \rho_i ) = -\rho_i \nabla \mathbf{u}_i, \ \ i=1,..N,\\ 
\end{equation}

\begin{equation}
\label{eq:motionSPHgas}
    \displaystyle\frac{ \partial \mathbf{v}}{\partial t}+ ( \mathbf{v} \cdot \nabla \mathbf{v} ) = -\displaystyle \frac{\nabla p}{\rho_{\rm g}} -\sum_i \varepsilon_i \displaystyle\frac{\mathbf{v}-\mathbf{u}_i}{t_{i}}+\mathbf{f}_{\rm g},
\end{equation}

\begin{equation}
\label{eq:motionSPHdust}
    \displaystyle\frac{ \partial \mathbf{u}_i}{\partial t}+ ( \mathbf{u}_i \cdot \nabla \mathbf{u}_i ) = \displaystyle\frac{\mathbf{v}-\mathbf{u}_i}{t_{i}}+\mathbf{f}_i, \ \ i=1,..N.\\
\end{equation}
Here $p$ is gas pressure, $\varepsilon_i=\displaystyle\frac{\rho_i}{\rho_{\rm g}}$ is fraction $i$ to gas mass ration, $\mathbf{f}_{\rm g}$, $\mathbf{f}_i$ are the acceleration due to forces acting on gas and dust fractions except for drag and pressure.

For the closure of the system (\ref{eq:contSPHgas})-(\ref{eq:motionSPHdust}), we need energy equation and equation of state for gas. In the paper, we consider two test problems: with complete equation for internal energy and with its particular case, i.e., isothermal approximation. 

We will use the complete equation for internal gas energy

\begin{equation}
\label{eq:innerEnSPH}
    \displaystyle \frac {\partial e}{\partial t}+ (\mathbf{v} \cdot \nabla e) = -\displaystyle \frac{p}{\rho_{\rm g}} \nabla \mathbf{v}-\sum_i \varepsilon_i \displaystyle\frac{(\mathbf{v}-\mathbf{u}_i)^2}{t_{i}}
\end{equation}
with the equation of state for ideal gas
\begin{equation}
\label{eq:eqidealState}
    p=\rho_{\rm g} e (\gamma-1),
\end{equation}
where $\gamma$ is the adiabatic exponent of gas. 

When solving the complete energy equation, we obtain system (\ref{eq:contSPHgas})-(\ref{eq:eqidealState}) consisting of ($2N+3$) partial differential equations (among them there is ($N+1$) vector equation, i.e., motion equations for gas and dust and ($N+2$) scalar equations, i.e., continuity equations for gas and dust and energy equation) and 1 algebraic equation. In the system, the unknown values are $\rho_{\rm g}$, $\rho_i$ ($N$ scalar values), $\mathbf{v}$ (the vector), $\mathbf{u}_i$ ($N$ vectors), $e$, $p$. 

For isothermal case, the analogs of (\ref{eq:innerEnSPH})-(\ref{eq:eqidealState}) will be 

\begin{equation}
\label{eq:econst}
    e=const,
\end{equation}
\begin{equation}
\label{eq:pisothermal}
    p=c_{\rm s}^2 \rho_{\rm g},
\end{equation}
where 
\begin{equation}
\label{eq:cspure}
c_{\rm s}=\sqrt{\displaystyle\frac{\partial p}{\partial \rho_{\rm g}}}=\sqrt{e (\gamma-1)}.    
\end{equation}

In the following section, we state two one-dimensional problems on the dynamics of gas-dust medium, which have reference solutions and are used to verify methods for calculation of momentum exchange \cite{LaibePrice2011Test,StoyanovskayaDustMultigrain,BenitezLlambay2019}.

\subsection{Test 1 DustyWave – propagation of acoustic waves}

The first test problem is the propagation of one-dimensional sound wave in isothermal gas-dust medium. This is a simple and demonstrative test that allows us to study how a numerical method reproduces a convective transfer simultaneously with the momentum exchange between the phases. The dynamics of isothermal gas-dust medium is described by the system (\ref{eq:contSPHgas})--(\ref{eq:motionSPHdust}) and (\ref{eq:econst})--(\ref{eq:pisothermal}) with $f_{\rm g}=0$, $f_i=0$, which is rewritten with variables $\mathbf{f}=(\rho_{g},\rho_i,v,u_i)$ as follows:

\begin{align}
    \frac{\partial \rho_{\rm g}}{\partial t}+\rho_{\rm g} \frac{\partial v}{\partial x} &= 0, \label{eq:DWS1} \\
    \frac{\partial \rho_i}{\partial t}+\rho_i \frac{\partial u_i}{\partial x} &= 0, \label{eq:DWS2} \\
    \frac{d v}{d t}+ \sum_{i=1}^{N}\frac{\rho_i}{\rho_{\rm g} t_i} (v-u_i)+\frac{c^2_s}{\rho_{\rm g}} \frac{\partial \rho_{\rm g}}{\partial x}&= 0, \label{eq:DWS3} \\
    \frac{d u_i}{d t}+\frac{1}{t_i}(u_i-v)&= 0. \label{eq:DWS4} 
\end{align}

In \ref{sec:type}, it is shown that that system (\ref{eq:DWS1})-(\ref{eq:DWS4}) is hyperbolic for any $\mathbf{f}$. Therefore, the Cauchy problem for it in the domain $x \in [0,L]$, $t \in [0,T]$ will have a unique solution. 

Let us state the Cauchy problem for the system (\ref{eq:DWS1})-(\ref{eq:DWS4}) 
\begin{equation}
\label{eq:initDustyWave}
    \mathbf{f}(x,t=0)=\mathbf{f}^0(x),
\end{equation}
\begin{equation}
\label{eq:boundaryDustyWave}
    \mathbf{f}(x=0,t)=\mathbf{f}(x=L,t),
\end{equation}
for $\mathbf{f}^0(x)$ to be continuous functions, $\rho^0(x) \neq const$. 

The constant solution
\begin{equation}
\label{eq:steadySolDW}
    \rho_{\rm g}(x,t)=\rho^0_{\rm g}, \quad \rho_{i}(x,t)=\rho^0_{i}, \quad v(x,t)=0, \quad u_i(x,t)=0
\end{equation}
satisfies (\ref{eq:DWS1})-(\ref{eq:DWS4}). Therefore, in a small neighborhood of (\ref{eq:steadySolDW}), the solutions of initial systems and the system linearized on (\ref{eq:steadySolDW}) are expected to be close. 

To study the properties of a numerical method, we will take the solution of linearized problem as a reference solution. This solution is found by the Fourier method. The overall procedure of obtaining the reference solution for dusty gas was described in \cite{LaibePrice2011} for monodisperse dust and in \cite{BenitezLlambay2019} for polydisperse dust. We used the prescription presented in \cite{BenitezLlambay2019} and developed SciLab program code that allows us to get the reference solution for specified parameters, initial and boundary conditions of the problem. The code is described in \ref{sec:dustyWaveMultigrain}, published freely\, and available at \verb|github.com| (see the link in subsection \ref{sec:intr_structure}). 

Note that the solution of linearized problem exists for all positive values of $t_i$. Therefore, the test makes it possible to study the properties of the new method both in a limiting case when for all fractions $t_i \ll 1$ and in the case of multi-scale stopping time values that is important in constructing AP schemes. 

\subsection{Test 2 DustyShock – Shock-tube problem}

The second test problem allows us to study the ability of a numerical method to reproduce the propagation of different types of waves: shock wave, rarefaction wave, and contact discontinuity in a polydisperse gas-dust medium. For the system (\ref{eq:contSPHgas})-(\ref{eq:motionSPHdust}) and (\ref{eq:innerEnSPH})-(\ref{eq:eqidealState}), we assume that the second term on the right-hand side of (\ref{eq:innerEnSPH}) equals 0 and $f_{\rm g}=0$, $f_i=0$. We formulate the Cauchy problem in the domain $x \in [0,L]$, $t \in [0,T]$. Let the mixture of gas and dust has zero velocity at the initial moment, but the values  $\mathbf{f}=(\rho_{g},\rho_i,e,p)$ have discontinuity
\begin{equation}
  f(x,t=0)|_{x<0.5L}=f^0_L, \quad f(x,t=0)|_{x>0.5L}=f^0_R.
\end{equation}
We take $T$ such that the waves propagating from the initial point of discontinuity  $x=0.5L$ do not reach the boundaries of the interval $x \in [0,L]$. Then on these boundaries, we give conditions for an unperturbed medium: 
\begin{equation}
\label{eq:boundaryDShock}
  f(x=0,t)=f^0_L, \quad f(x=L,t)=f^0_R \quad v(x=0,t)=v(x=L,t)=0, \quad u_i(x=0,t)=u_i(x=L,t)=0.
\end{equation}
  
The analytical solution of such a problem for pure gas is known in the whole domain of parameters. This solution is widely used to test the numerical methods for gas-dynamic equations \cite{Sod1978}. For gas with monodisperse dust particles in \cite{LaibePrice2011Test}, it is shown that in the extreme case of short stopping time, gas-dust medium behaves as a heavier gas with efficient reduced sound speed 
\begin{equation}
    c^*_{\rm s}=\displaystyle\frac{c_{\rm s}}{\sqrt{1+\displaystyle\frac{\rho_{\rm d}}{\rho_{\rm g}}}},
\end{equation}
where $\rho_{\rm d}$ is the density of monodisperse dust. Therefore, the reference solution is constructed by replacing $c_{\rm s}$ by $c^*_{\rm s}$ in the solution for pure gas. In \ref{sec:analysis}, we demonstrate that in the extreme case of short stopping time values for the medium with polydisperse particles 
\begin{equation}
    c^*_{\rm s}=\displaystyle\frac{c_{\rm s}}{\sqrt{1+\sum_i \varepsilon_i}}.
\end{equation}  
Therefore, to obtain the reference solution of the problem DustyShock with polydisperse dust, we can use a program for monodisperse dust supposing that $\rho_{\rm d}=\sum_i \rho_i$. In our paper, we applied the code SPLASH published by the authors~\cite{SPLASH}.

%%%%%%%%%%%%%%%%%%%%%%%%%%%%%%%%%%%%%%%%%%%%%%%%%%%%%%%%%%%%%%%%%

\section{Numerical algorithm}
\label{sec:NumMethod}

In Section \ref{sec:MethodIdea}, we present the ideas forming the basis of the proposed method to solve motion equations (\ref{eq:motionSPHgas})-(\ref{eq:motionSPHdust}). Section \ref{sec:methodformulas} contains the calculation formulas for solving the complete system of equations (\ref{eq:contSPHgas})-(\ref{eq:innerEnSPH}) constructed on the base of these ideas.

\subsection{The main ideas of drag computing}
\label{sec:MethodIdea}
We suppose that the problem solved on the base of equations (\ref{eq:contSPHgas})-(\ref{eq:innerEnSPH}) is dynamical, hence, it is reasonable to apply time-explicit standard SPH-approximation \cite{SPH} of equations of continuity and energy and all forces except for drag in the motion equations. 

\begin{figure}[!t]
\begin{center}
\includegraphics{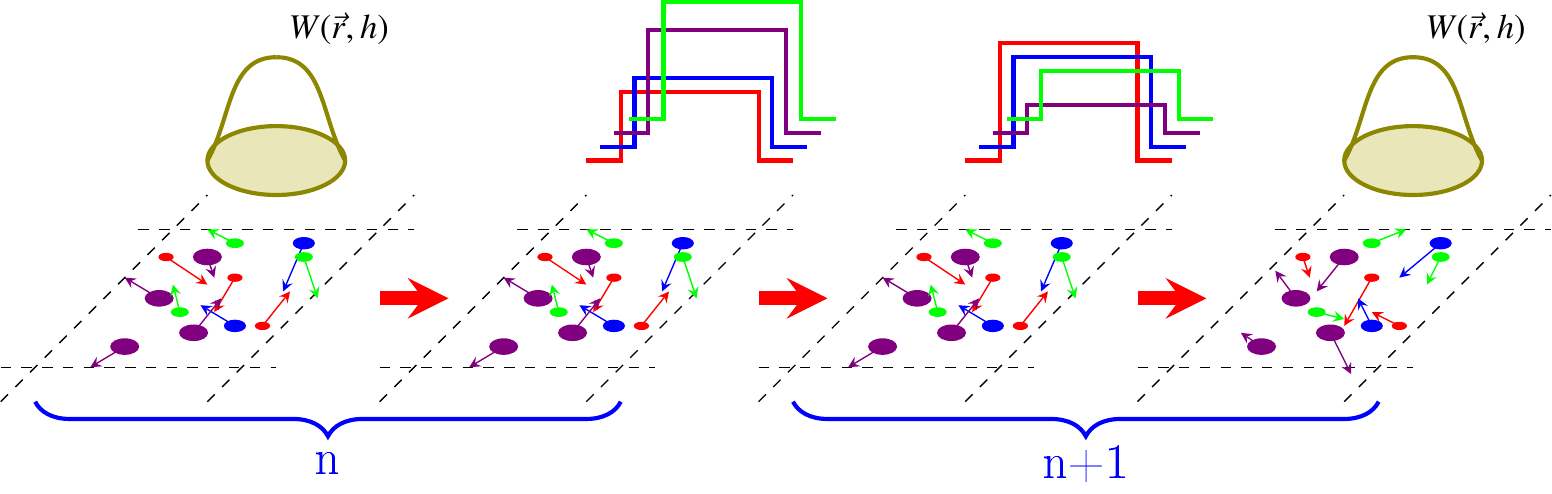}
\caption{Calculation of velocity in the method IDIC for two-phase polydisperse medium. The example with 3 dust fractions is shown. Step 1 is the calculation of average velocity values in cells at a current time step. Step 2 is the calculation of the velocity values in cells at the next time step by solving the system of linear algebraic equations in $O(N^2)$ arithmetic operations. Step 3 is the calculation of velocity values of the particles at the next time step.} 
\label{fig:IDIC}
\end{center}
\end{figure}

Here, to calculate relative velocity values between gas and particle of different fractions, we will use the method of cell-averaging, which is efficient in the case of monodisperse dust \cite{IDIC2018}. In \cite{IDIC2018}, from formulas (44)-(45) and (10), we concluded that this method preserves the momentum in a cell with machine precision. The advantage of cell-averaging approach is that it is easily extended from monodisperse dust to polydisperse one. The main steps of calculation are illustrated in Fig.~\ref{fig:IDIC}. By analogy with approach \cite{IDIC2018}, we cover the calculation domain with a celled grid of arbitrary shape and calculate gas velocity values $\mathbf{v}_*$ and dust fractions $\mathbf{u}_{i*}$ in each grid cell as arithmetic mean of velocity of all model particles of this type. We write the motion equations in such a way that in the calculation of drag acting on the particle, a cell-averaged velocity is involved and vice versa:

\begin{equation}
\label{eq:motionAVEgas}
    \displaystyle\frac{ \partial \mathbf{v}}{\partial t}+ ( \mathbf{v} \cdot \nabla \mathbf{v} ) = -\displaystyle \frac{\nabla p}{\rho_{\rm g}} -\sum_i \varepsilon_i \displaystyle\frac{\mathbf{v}-\mathbf{u}_{i*}}{t_{i}}+\mathbf{f}_{\rm g},
\end{equation}

\begin{equation}
\label{eq:motionAVEdust}
    \displaystyle\frac{ \partial \mathbf{u}_i}{\partial t}+ ( \mathbf{u}_i \cdot \nabla \mathbf{u}_i ) = \displaystyle\frac{\mathbf{v_*}-\mathbf{u}_i}{t_{i}}+\mathbf{f}_{i}, \ \ i=1,..N.\\
\end{equation}

To find a solution with time step, which is defined by not small stopping time $t_i$, but the Courant condition, we approximate implicitly relaxation terms in system (\ref{eq:motionAVEgas})-(\ref{eq:motionAVEdust})as follows: 

\begin{equation}
\label{eq:motionAVEIMPLICITgas}
    \displaystyle\frac{ \mathbf{v}^{n+1}-\mathbf{v}^{n}}{\tau}+ ( \mathbf{v}^n \cdot \nabla \mathbf{v}^n ) = -\displaystyle \frac{\nabla p^n}{\rho^n_{\rm g}} -\sum_i \varepsilon^n_i \displaystyle\frac{\mathbf{v}^{n+1}-\mathbf{u}^{n+1}_{i*}}{t_{i}}+\mathbf{f}^n_{\rm g},
\end{equation}

\begin{equation}
\label{eq:motionAVEIMPLICITdust}
    \displaystyle\frac{ \mathbf{u}^{n+1}_i-\mathbf{u}^n_i}{\tau}+ (\mathbf{u}^n_i \cdot \nabla \mathbf{u}^n_i ) = \displaystyle\frac{\mathbf{v}^{n+1}_*-\mathbf{u}^{n+1}_i}{t_{i}}+\mathbf{f}^n_{i}, \ \ i=1,..N.\\
\end{equation}

 From equations (\ref{eq:motionAVEIMPLICITgas})-(\ref{eq:motionAVEIMPLICITdust}), we exclude the terms on the left-hand side responsible for the transfer and we write the obtained equations for each particle from the cell. Next, we sum up them for the particles of each type and divide them by the number of particles of such type. Assume that
\begin{equation}
    \mathbf{A}_{\rm g}^n=-\displaystyle \left(\frac{\nabla p^n}{\rho^n_{\rm g}}\right)+\mathbf{f}^n_{\rm g}, \quad \mathbf{A}_i^n=\mathbf{f}^n_{i}.
\end{equation}

We obtain the system of equations to define the cell-averaged velocity values at the step ($n+1$):

\begin{equation}
\label{eq:motionCELLgas}
    \displaystyle\frac{ \mathbf{v}_*^{n+1}-\mathbf{v}_*^{n}}{\tau} = \mathbf{A}^n_{\rm g*} -\sum_i \varepsilon_{*i}^n \displaystyle\frac{\mathbf{v}_*^{n+1}-\mathbf{u}^{n+1}_{i*}}{t_{i}},
\end{equation}

\begin{equation}
\label{eq:motionCELLdust}
    \displaystyle\frac{ \mathbf{u}^{n+1}_{*i}-\mathbf{u}_{*i}^n}{\tau} = \mathbf{A}^n_{i*}+\displaystyle\frac{\mathbf{v}^{n+1}_*-\mathbf{u}^{n+1}_{*i}}{t_{i}}, \ \ i=1,..N.\\
\end{equation}

For small $t_i$, the relaxation terms are a small difference between two large values $\mathbf{v}_*$ and $\mathbf{u}_{i*}$, which is divided by a small value $t_i$. In floating-point arithmetic, the direct calculation of such a term may cause the loss of accuracy. Therefore, we replace the system (\ref{eq:motionCELLgas})-(\ref{eq:motionCELLdust}) by the equivalent one written in variables of barycentric and relative velocity values as
\begin{equation}
\label{eq:vuToxy}
    \mathbf{w}=\mathbf{v}_*+\sum_i \varepsilon_i \mathbf{u}_{*i}, \quad \mathbf{w}_i=\mathbf{v}_*-\mathbf{u}_{*i},
\end{equation}

\begin{equation}
\label{eq:motionY}
    \displaystyle\frac{ \mathbf{w}^{n+1}-\mathbf{w}^{n}}{\tau} = \mathbf{A}^n_{\rm g*}+\sum_i \varepsilon_i^*\mathbf{A}^n_{i*},
\end{equation}

\begin{equation}
\label{eq:motionX}
    \displaystyle\frac{ \mathbf{w}_i^{n+1}-\mathbf{w}^n_i}{\tau} = \mathbf{A}^n_{\rm g*}-\mathbf{A}^n_{i*} -\displaystyle\frac{\varepsilon_{*i}+1}{t_{i}}\mathbf{w^{\textcolor{red}{n+1}}}_i-\sum_{j \neq i} \frac{\varepsilon_{*j}}{t_j}\mathbf{w^{\textcolor{red}{n+1}}}_j, \ \ i=1,..N.\\
\end{equation}

Equation (\ref{eq:motionY}) can be solved apart from the others, but (\ref{eq:motionX}) in one-dimensional case is the linear system of dimension $N \times N$ with a filled matrix of special structure. The substitution
\begin{equation}
\label{eq:xtoz}
z_i^n=\displaystyle\frac{\tau \varepsilon_i}{t_i}w_i^n, \quad 
b_i=\displaystyle\frac{t_i+\tau}{\varepsilon_i\tau}
\end{equation}
yields the system $B\mathbf{z}^{n+1}=\tau(\mathbf{A}^n_{\rm g*}-\mathbf{A}^n_{i*})+\displaystyle\frac{t_i}{\tau \varepsilon_i}\mathbf{z}_i^n$ with matrix $B$ that is easily inverted analytically: 
\begin{equation}
\label{eq:ntoSolve}
B=
\begin{pmatrix}
1+b_1 & 1 & ... & 1 \\
1 & 1+b_2 & ... & 1 \\
... & ... & ... &... \\
1 &  1    & ... & 1+b_N
\end{pmatrix},
\quad
B^{-1}=-\displaystyle\frac{1}{\beta}
\begin{pmatrix}
\displaystyle\frac{1-b_1 \beta}{b_1^2} & \displaystyle\frac{1}{b_1 b_2} & ... & \displaystyle\frac{1}{b_1 b_N} \\
\displaystyle\frac{1}{b_2 b_1} & \displaystyle\frac{1-b_2 \beta}{b_2^2} & ... & \displaystyle\frac{1}{b_2 b_N} \\
... & ... & ... &... \\
\displaystyle\frac{1}{b_N b_1} &  \displaystyle\frac{1}{b_N b_2}    & ... & 
\displaystyle\frac{1-b_N \beta}{b_N^2}
\end{pmatrix},
\quad
\beta=1+\displaystyle\sum_{i=1}^N\displaystyle\frac{1}{b_i}.
\end{equation}
Thus the solution of system (\ref{eq:motionX}) can be found in $O(N^2)$ arithmetic operations if $\varepsilon_i>0 \quad \forall i$. 
  
Note that the similar problem of finding the velocity values at the next time moment from the solution of the system of linear algebraic equations arises in constructing grid methods for a two-phase polydisperse medium. The authors \cite{BenitezLlambay2019,Li2017} solve the system of equation (\ref{eq:motionCELLgas})-(\ref{eq:motionCELLdust}) by using Gaussian elimination. In the Gaussian method, $O((N+1)^3)$ arithmetic operations are required for solving the system of dimension $(N+1)\times(N+1)$. The use of quadratic complexity algorithm, which we proposed, will essentially reduce the computational cost in those applications where a large number of fractions is considered (for example, the authors \cite{Drajkowska2019} use 150 dust fractions).  

After finding $\mathbf{w}^{n+1}$, $\mathbf{w}_i^{n+1}$, we define the velocity values of gas and dust fractions from (\ref{eq:vuToxy}).
According to the velocity values of gas and dust found in each cell, we restore the velocity values of individual SPH-particles by equations (\ref{eq:motionAVEIMPLICITgas})-(\ref{eq:motionAVEIMPLICITdust}) recorded for individual particles.

Let us estimate the number of operations required for drag computing in standard SPH 'particle-particle' way and suggested drag-in-cell or 'particle-mesh' way. Let $N_{\rm cell}$ be the number of grid cells, $N_{\rm neib}$ be the typical number of neighbors of the same sort (gas or $i$-th fraction of dust) inside smoothing kernel, $N_{\rm SPH}$ be the number of SPH particles of every sort. For simplicity, let us start with explicit drag approximation. To extend the Monaghan and Kocharyan way to polydisperse case, one needs to find neighbors for all particles of every sort using $O_1=O((N+1)\times (N_{\rm SPH} \ln{N_{\rm SPH}}))$ operations, and then for every particle to sum up contribution of neighbors of all sorts using $O_2=O((N+1) \times N_{\rm SPH} \times N_{\rm neib} \times (N+1))$ operations. To implement a new drag-in-cell way, one needs to fill grid cells with averaged velocities using $O_3=O\left(N_{\rm SPH} \times (N+1)+N_{\rm cell} \right)$ operations (Step 1 of Fig.\ref{fig:IDIC}), then to update velocity values in cells with $O_4=O(N_{\rm cell} \times (N+1))$ operations (Step 2 of Fig.\ref{fig:IDIC}) and then to update velocity values of individual particles with $O_5=O(N_{\rm SPH} \times (N+1))$ operations (Step 3 of Fig.\ref{fig:IDIC}). To make the drag approximation implicit in the Monaghan and Kocharyan way, one needs to make iterations repeating $O_2$ step several times until convergence. To make the drag approximation implicit in a drag-in-cell way, one needs only to increase the complexity of $O_4$ to $O(N_{\rm cell} \times (N+1)^2)$. It is clearly seen that the complexity of $O_4$, even in case of implicit approximation, is less than complexity of $O_2$ with explicit approximation under natural condition $N_{\rm cell}<N_{\rm SPH} \times N_{\rm neib}$. This efficiency is achieved owing to treating a group of particles inside one cell as a superparticle with single velocity and using of extra memory to store grid values.

\subsection{Calculation formulas}
\label{sec:methodformulas}

\begin{figure}[!t]
\includegraphics{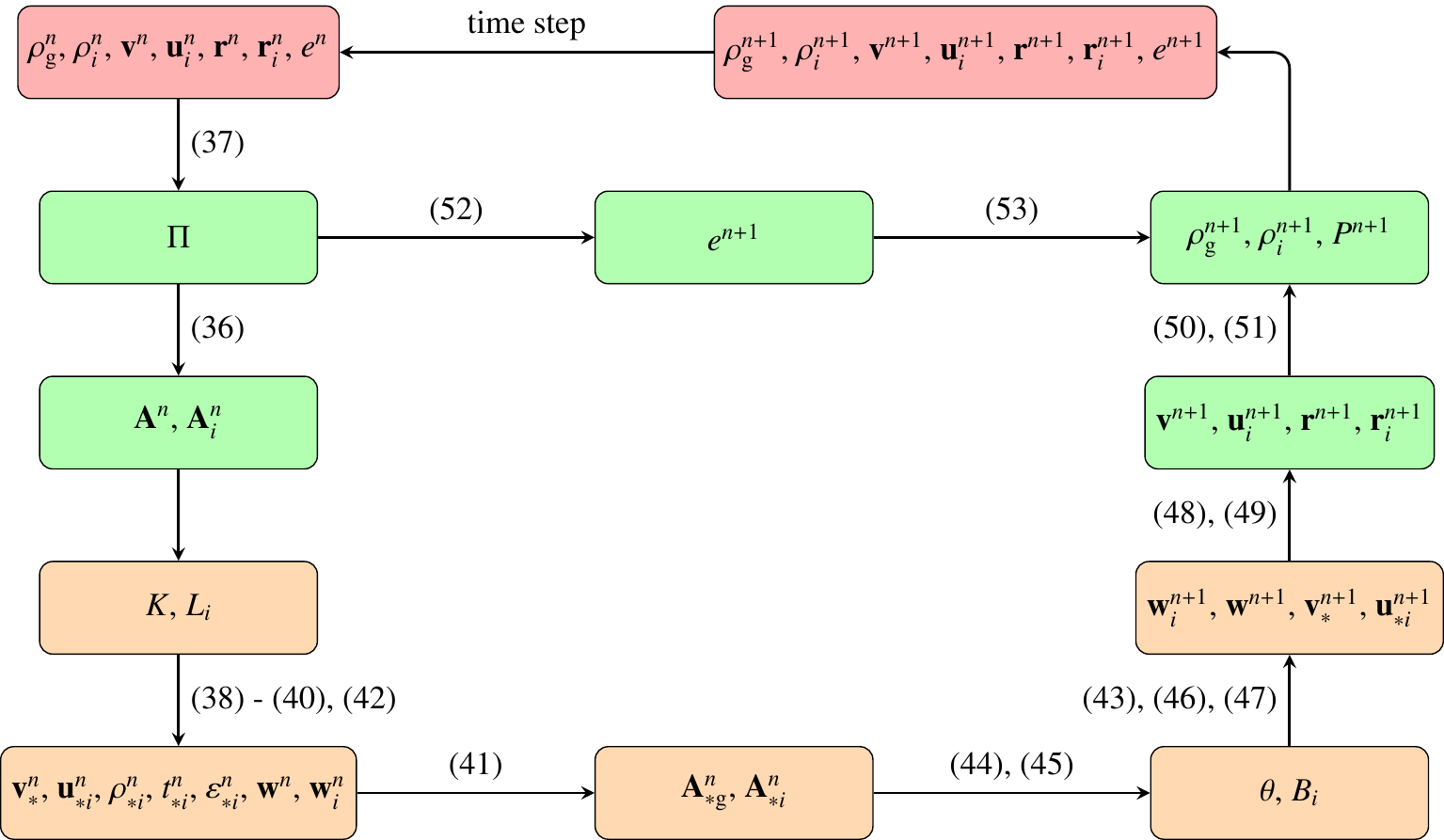}
\caption{The flowchart of one time step.} \label{fig:Scheme}
\end{figure}

Let the coordinates, velocity, density of gas and dust particles and the internal energy of gas particles be known at time moment $n$. We will find these values at time moment $(n+1)$ calculating the velocity values of the particles by using a new hybrid method. To calculate coordinates, density, entropy and pressure, standard SPH-approximation is applied. Fig.~\ref{fig:Scheme} illustrates the flowchart of one time step. Input and output data are shown in red, the values calculated in particles are shown in green, and the values calculated in cells are shown in orange. 

We use $a,b$ as subscripts for gas particles and $m,l$ are subscripts for dust particles. Thus, for example $\rho_{a, \rm g}, \mathbf{v}_a, \mathbf{r}_a, P_a$ will denote density, velocity, radius-vector, gas-kinetic  pressure in $a$-numbered  particle, respectively; $\rho_{il}, \mathbf{u}_{il}, \mathbf{r}_{il}$ are density, velocity and radius-vector of $i$-th dust fraction in $l$-numbered particle.

\subsubsection{Calculation of velocity values of gas and dust particles by using the hybrid method "particle-particle" and "particle-cell".}
\label{sec:IDICMultigrain}

STEP 0. At each time moment, we divide the whole computational domain into non-overlapping volumes in such a way that the union of these volumes coincides with the whole domain. Let $K$ gas particles of the same mass $m_{\rm g}$ and $L_i$ dust particles of fraction $i$ of the same mass $m_i$, here $K>0$, $L_i>0$ be in a separate volume. 

STEP 1. We calculate the acceleration due to all forces except for drag, which act on the particles at time moment $n$:

\begin{equation}
\label{eq:Aparticle}
\mathbf{A}^{n}_a  = - \displaystyle \sum_b m_g\left(\frac{P^{n}_b}{(\rho^n_{b,\rm{g}})^2} + \frac{P^{n}_a}{(\rho^n_{a,\rm g})^2} + \Pi_{ab} \right) \bigtriangledown_a W_{ab}, \quad  \mathbf{A}_{ij}^n=0,   
\end{equation}
where
\begin{equation}
\label{eq:Visc}
\Pi_{ab} = 
\begin{cases}
\displaystyle \frac{-\alpha c_{ab} \mu_{ab} + \beta \mu^2_{ab}}{\rho_{ab}}, &\text{if $ \mathbf{v}_{ab} \mathbf{r}_{ab} < 0$,}\\
0, &\text{if $ \mathbf{v}_{ab} \mathbf{r}_{ab} > 0$,}
\end{cases}
\end{equation}

$\mu_{ab} = \displaystyle \frac{h \mathbf{v}_{ab} \mathbf{r}_{ab}}{\mathbf{r}^2_{ab} + \nu^2}$,
 $\mathbf{v}_{ab} = \mathbf{v}_a - \mathbf{v}_b$, \quad $\mathbf{r}_{ab} = \mathbf{r}_a - \mathbf{r}_b$, \quad $\rho_{ab} = \displaystyle \frac{1}{2} (\rho_{a,\rm g} + \rho_{b,\rm g})$, \quad $c_{ab} = \displaystyle \frac{1}{2} (c_a + c_b)$ \quad and \quad $c_a = \sqrt[]{\displaystyle \frac{\gamma P_a}{\rho_{a,\rm g}}}$, $\nu$ is a clipping constant or limiter, which does not allow a denominator to become 0 for approaching SPH particles. 

STEP 2. We calculate cell-averaged values at time moment $n$:

\begin{equation}
\label{eq:VUvolaverage}
\mathbf{v}_*^n=\displaystyle\frac{\sum_{a=1}^K \mathbf{v}_a^n}{K}, \quad \mathbf{u}_{*i}^n=\displaystyle\frac{\sum_{l=1}^{L_i} \mathbf{u}^n_{il}}{L_i},
\end{equation}

\begin{equation}
t_{*i}^n=\displaystyle\frac{\sum_{l=1}^{L_i} t^n_{il}}{L_i},   
\end{equation} 

\begin{equation}
\label{eq:aveEpsilon}
\varepsilon^n_{*i}=\displaystyle\frac{m_i L_i}{m_{\rm g} K}, 
\end{equation} 

\begin{equation}
\label{eq:Avolaverage}
\mathbf{A}_{*g}^n=\displaystyle\frac{\sum_{a=1}^K \mathbf{A}_a^n}{K}, \quad \mathbf{A}_{*i}^n=\displaystyle\frac{\sum_{l=1}^{L_i} \mathbf{A}^n_{il}}{L_i},
\end{equation}

STEP 3. We find the velocity values averaged over each cell at the next time moment $\mathbf{u}_{*i}^{n+1}$, $\mathbf{v}_*^{n+1}$.
 
For this purpose, we assume that
\begin{equation}
\label{eq:newNxy}
\mathbf{w}^n=\mathbf{v}_*^n+\displaystyle \sum_i \varepsilon_{*i}^n \mathbf{u}_{*i}^n, \quad \mathbf{w}_i^n=\mathbf{v}_*^n-\mathbf{u}_{*i}^n.
\end{equation}

Then 
\begin{equation}
\label{eq:yN}
\mathbf{w}^{n+1}=\mathbf{w}^n+\tau(\mathbf{A}_{*g}^n+\sum_i \varepsilon_{*i}^n \mathbf{A}_{*i}^n),
\end{equation}

\begin{equation}
\label{eq:Bi}
    B_i=\displaystyle\frac{t_{*i}^n+\tau}{\varepsilon_{*i}^n\tau},
\end{equation}

\begin{equation}
\label{eq:theta}
    \theta=1+\displaystyle\frac{1}{B_1}+\frac{1}{B_2}+...\frac{1}{B_N},
\end{equation}

\begin{equation}
\label{eq:Nsolution}
\mathbf{w}^{n+1}_i=-\displaystyle\frac{t_{*i}^n}{\tau \varepsilon_{*i}^n \theta}
\left[ \frac{(1-B_i \theta)(\mathbf{w}_i^n+\tau(\mathbf{A}_{*g}^n-\mathbf{A}_{*i}^n))}{B_i^2}+
\sum_{j\ne i}\frac{\mathbf{w}_j^n+\tau(\mathbf{A}_{*g}^n-\mathbf{A}_{*j}^n)}{B_i B_j} \right],
\end{equation}

\begin{equation}
\label{eq:vufinal}
\mathbf{v}_*^{n+1}=\displaystyle\frac{\mathbf{w}^{n+1}+\sum_i \varepsilon_{*i}^n \mathbf{w}^{n+1}_i}{1+\sum_i \varepsilon_{*i}^n}, \ \ \mathbf{u}_{i*}^{n+1}=\displaystyle\frac{\mathbf{w}^{n+1}-(1+\sum_{j \neq i}\varepsilon_{*j}^n)\mathbf{w}^{n+1}_i+\sum_{j \neq i} \varepsilon_{*j}^n \mathbf{w}^{n+1}_j}{1+\sum_j \varepsilon_{*j}^n}.
\end{equation}

STEP 4. We calculate the velocity values of each gas and dust particle at time moment $n+1$.

\begin{equation}
\label{eq:VUNew}
\left(\frac{1}{\tau}+\sum_{i=1}^N \frac{\varepsilon^n_{*i}}{t^n_ {*i}}\right)\mathbf{v}^{n+1}_a =\frac{\mathbf{v}^n_a}{\tau}+ \sum_{i=1}^N  \frac{\varepsilon^n_{*i}}{t^n_{*i}} \mathbf{u}^{n+1}_{*i} + \mathbf{A}_a^n,  \quad 
\left(\frac{1}{\tau}+\frac{1}{t^n_{*i}}\right)\mathbf{u}^{n+1}_{il} =\frac{\mathbf{u}^n_{il}}{\tau}+ \frac{1}{t^n_{*i}} \mathbf{v}^{n+1}_*+\mathbf{A}_{ij}^n.
\end{equation}

\subsubsection{Calculation of radius-vectors, density, internal energy and pressure by using SPH method}

We calculate new radius-vectors using velocity values $\mathbf{v}^{n+1}$ 

\begin{equation}
\label{eq:positionUPD}
\mathbf{r}_{a, \rm g}^{n+1}=\mathbf{r}_{a, \rm g}^n+\tau \mathbf{v}_a^{n+1}, \quad \mathbf{r}_{il}^{n+1}=\mathbf{r}_{il}^n+\tau \mathbf{u}_{il}^{n+1}.   
\end{equation}

The density of dust and gas is calculated not by approximation of continuity equations, but by interpolation of point masses located in particles
\begin{equation}
\label{eq:contin_gas}
\rho^{n+1}_{a, \rm g} = m_{\rm g} \sum_b W^{n+1}_{ab},
\end{equation}
\begin{equation}
\label{eq:contin_dust}
\rho^{n+1}_{il} = m_{i} \sum_m W^{n+1}_{lm}.
\end{equation}

We approximate Eq.~(\ref{eq:innerEnSPH}) assuming that the second term on its right-hand side equals zero that yields
\begin{equation}
\label{eq:energySPH}
e_a^{n+1} = e_a^n+\tau \left(\displaystyle\frac{m_{\rm g} P_a}{(\rho^n_{a, \rm g})^2} \sum_b (\mathbf{v}^n_a - \mathbf{v}^n_b) \nabla_a W^n_{ab} + \frac{m_{\rm g}}{2} \sum_b \Pi_{ab} (\mathbf{v}^n_a - \mathbf{v}^n_b) \nabla_a W^n_{ab}\right),
\end{equation}
and then with allowance for internal energy found, we calculate pressure as
\begin{equation}
\label{eq:eqidealStateSPH}
    P_a^{n+1}=\rho_{a,\rm g}^{n+1} e_a^{n+1} (\gamma-1).
\end{equation}

\section{Calculation results}
\label{sec:results}

In this section, we study the properties of the proposed method IDIC to calculate drag between gas and dust. The focus of our attention is to define if the proposed method of calculation of interphase interaction allows one to simulate the dynamics of two-phase polydisperse medium using the same spatial and time resolution as the dynamics of pure gas.  

Test calculations are the solutions to DustyWave and DustyShock problems. In \ref{sec:setup}, we define the initial data, physical and numerical parameters of test problems. 

It is known from the simulation of dynamics of medium with monodisperse inclusions that short stopping time requires high numerical resolution. In \ref{sec:smalltstop}, we study the models, in which all dust fractions have small $t_i$. When passing from monodisperse approximation to polydisperse one, the availability to consider multi-scale stopping time values of dust fraction is fundamentally new.  In \ref{sec:multiscaletstop}, we show the calculation results for $t_i$, which differ by the orders of magnitude; moreover, one of dust fractions moves almost at gas velocity and the velocity of others differs from that of gas.   

\subsection{Initial data, physical and calculation parameters of test problems}
\label{sec:setup}

Numerical artefacts become more clearly defined at high dust to gas mass ratio. Therefore, in all calculations, we will specify the total dust density equal to gas density. Moreover, in all calculations, $N_{\rm SPH}$, which is the number of particles modeling gas, was equal to the number of particles modeling each dust fraction, so the total number of particles is $(N+1) \times N_{\rm SPH}$. 

All calculations in the paper were performed by using constant smoothing radius $h$ and constant time step $\tau$, which satisfies the Courant condition for pure gas
\begin{equation}
\label{eq:courant}
\tau \le \displaystyle \frac{h \cdot {\rm CFL}}{\max(c_{\rm s}, u_i, v)},
\end{equation}
where ${\rm CFL}$ is the Courant parameter. In most cases ${\rm CFL}$ equal to 0.5 or 0.25 were set. To calculate the velocity, we used a motionless uniform grid of cell size  $h_{\rm cell}$ equal to $0.5h$. 

For each calculation, the number of model particles $N_{\rm SPH}$, smoothing radius $h$ and the size of cell $h_{\rm cell}$, and time step $\tau$ are given in Table~\ref{tab:runsetup}.

We used a standard cubic spline as a kernel for the one-dimensional case 
\begin{equation}
\label{eq:kernel}
W^{n}_{ab} = W(|r^{n}_a - r^{n}_b|, h) = W^{n}(q) = \frac{2}{3 h}
\begin{cases}
  1 - \frac{3}{2} q^2 + \frac{3}{4} q^3, &\text{if $ \displaystyle 0 \le q \le 1$,}\\
  \displaystyle \frac{1}{4}(2 - q)^3, &\text{if $\displaystyle 1 \le q \le 2$,}\\
  0, &\text{otherwise;}
\end{cases}
\end{equation}
where $q = \displaystyle \frac{|r^{n}_a - r^{n}_b|}{h}$. 

\begin{table*}
\begin{minipage}{160mm}
\caption{Physical and numerical parameters of DustyWave (DW1 - DW3) and DustyShock (DS1-DS9) problems.}
\centering
 \label{tab:runsetup}
 \begin{tabular}{cccccccc}
 \hline
    Name & $N$ & $t_i$ or $s_i$ & $\varepsilon_i$ & $N_{\rm SPH}$ & $h$ & $\tau$ & Method \\
  \hline
  DW1 & 3 & $t_i=0.1,0.2,0.4$ & $\varepsilon_i=0.33,0.33,0.33$ & 600 & 0.01 & $5 \times 10^{-3}$ & IDIC \\
 
  DW2 & 3 & $t_i=10^{-2},10^{-3},10^{-4}$ & $\varepsilon_i=0.33,0.33,0.33$ & 600 & 0.01 & $5 \times 10^{-3}$ & IDIC \\
  
  \textcolor{magenta}{DW3} & 3 & $t_i=10^{-2},10^{-3},10^{-4}$ & $\varepsilon_i=0.33,0.33,0.33$ & \textcolor{magenta}{30} & \textcolor{magenta}{0.1} & $5 \times 10^{-3}$ & \textcolor{magenta}{IDIC} \\
  
  \hline
  DS1 & 1 & $s_i=10^{-4}$ & $\varepsilon_i=1$ & 2100 & 0.01 & $5 \times 10^{-3}$ & IDIC \\
  DS2 & 1 & $s_i=10^{-4}$ & $\varepsilon_i=1$ & 2100 & 0.01 & $5 \times 10^{-5}$ & MK \\
  DS3 & 1 & $s_i=10^{-4}$ & $\varepsilon_i=1$ & 21000 & 0.001 & $5 \times 10^{-4}$ & IDIC \\
  DS4 & 1 & $s_i=10^{-4}$ & $\varepsilon_i=1$ & 21000 & 0.001 & $5 \times 10^{-5}$ & MK \\
  
  \hline
  DS5 & 2 & $s_i=10^{-3},10^{-4}$ & $\varepsilon_i=0.01,0.99$ & 1180 & 0.01 & $5 \times 10^{-3}$ & IDIC \\
  
  DS6 & 2 & $s_i=10^{-3},10^{-4}$ & $\varepsilon_i=0.5,0.5$ & 1180 & 0.01 & $5 \times 10^{-3}$ & IDIC \\
  
  \hline
  DS7 & 3 & $s_i=10^{-3},10^{-2},10^{-1}$ & $\varepsilon_i=0.33,0.33,0.33$ & 1180 & 0.02 & $5 \times 10^{-3}$ & IDIC \\
  
  DS8 & 3 & $s_i=10^{-3},10^{-2},10^{-1}$ & $\varepsilon_i=0.33,0.33,0.33$ & 2360 & 0.01 & $2.5 \times 10^{-3}$ & IDIC \\
  
  DS9 & 3 & $s_i=10^{-3},10^{-2},10^{-1}$ & $\varepsilon_i=0.33,0.33,0.33$ & 7086 & 0.005 & $1.25 \times 10^{-3}$ & IDIC \\
 
\end{tabular}
\end{minipage}
\end{table*}

\subsubsection{DustyWave}

We publish the results of 3 runs DW1, DW2 and DW3 of DustyWave problem with 3 dust fractions. In the system (\ref{eq:DWS1})-(\ref{eq:DWS4}), values $t_i$, $c_{\rm s}$ are parameters. We specify $c_{\rm s}=1$ identical for all cases, the values $t_i$ are different and presented in Table~\ref{tab:runsetup}. 

We formulate the problem at interval $[0,1]$ assuming that initial values of density and velocity are the sum of background value and perturbations of small-amplitude  
\begin{equation}
 \rho_{\rm g}(x)=\rho_{\rm g}^0+\delta\rho_{\rm g}(x), \quad \rho_i(x)=\rho^0_i+\delta\rho_j(x), \quad v(x)=v^0+\delta v(x), \quad u_i(x)=u_i^0+\delta u_i(x).
\end{equation}
In all calculations, we will specify $v^0=u_i^0=0$, $\rho_{\rm g}^0=1$ and initial values $\varepsilon_i$ (are shown in Table~\ref{tab:runsetup}) using them in defining $\rho^0_i$. We will give perturbations as $\delta(x)=A(\varphi \cos{(2 \pi k x)}+\chi \sin{(2 \pi k x)})$, where amplitude $A$ and wave number $k$ are the parameters identical for all values of density and velocity, $\varphi$, $\chi$ are the coefficients defined for each variable by using the code we publish as supporting material. For convenient reproducing of the results, we present initial perturbations of density and velocity values for DW1, DW2, and DW3 in Table ~\ref{tab:DustyWaveInit}. 

In order to obtain density equal to $\rho_{\rm g}(x)=\rho^0_{\rm g}+\delta\rho_{\rm g}(x)$ on interval $[0,1]$, we used a recurrent procedure of locating model particles. The first particle was placed at the origin of coordinates $x_1=0$, the coordinate of the next one was determined from the relation   
$$\displaystyle \int\limits_{x_i}^{x_i + \Delta x_i} \rho_{\rm g}(x) dx = \frac{1}{N_{\rm SPH}}\displaystyle \int\limits_{0}^{1} \rho_{\rm g}(x) dx .$$
After finding all $x_i$, each particle was shifted right by a value $\Delta x_i$. The particles of dust fractions were placed in analogous way.

We are interested in the solution of the problem in the domain $x \in [0,1]$, $t \in [0,2]$. To obtain it, the calculation was performed in the extended domain of space $x \in [-4,4]$. For this purpose, at the initial time moment, gas and dust particles placed on $[0,1]$ were copied on $[-1,0]$, $[1,2]$, etc. Such a location of particles resulted in continuous periodic distributions of density values of gas and dust fractions. Further, the equations were solved for all particles in the domain $[-4,4]$. Such a calculation in the extended domain automatically guaranteed fulfillment of periodical boundary conditions (\ref{eq:boundaryDustyWave}) for $[0,1]$. 

Calculations of DW1, DW2, and DW3 were performed without artificial viscosity. 

\begin{table*}
\begin{minipage}{160mm}
\caption{Initial perturbations of background values of density $\rho^0_{\rm g}=1$, $\rho^0_{i}=0.3333$ and velocity $v^0=u^0_i=0$ in DustyWave models DW1, DW2 and \textcolor{magenta}{DW3}  with $A=0.0001$, $c_{\rm s}=1$, $k=1$, $x \in [0,1] $.}
\centering
 \label{tab:DustyWaveInit}
 \begin{tabular}{ccc}
 \hline
    & $t_1 = 10^{-2}$, $t_2 = 10^{-3}$, $t_3 = 10^{-4}$ (DW2, \textcolor{magenta}{DW3}) &  $t_1 = 0.1$, $t_2 = 0.2$, $t_3 = 0.3$ (DW1) \\
  \hline
  $\rho_{\rm g}$      & $A \cos(2\pi x)$   & $A \cos(2\pi x)$   \\
  $v$ & $A(-0.707212 \cos(2\pi x) + 0.0029033 \sin(2\pi x))$  &   $A(-0.7852741 \cos(2\pi x) + 0.1267991 \sin(2\pi x))$  \\
\hline
 $\rho_1$ & $A(0.3327036 \cos(2\pi x) + 0.0147865 \sin(2\pi x))$  &  $A (0.2813014 \cos(2\pi x) + 0.1508098 \sin(2\pi x))$ \\
 $\rho_2$ & $A(0.3332995\cos(2\pi x) + 0.0014811\sin(2\pi x))$  &  $A (0.1667321 \cos(2\pi x) + 0.1957177 \sin(2\pi x))$ \\     
 $\rho_3$ & $A(0.3333005\cos(2\pi x) + 0.0001481\sin(2\pi x))$  &  $A (0.0520914 \cos(2\pi x) + 0.1508957 \sin(2\pi x)) $
\\
 $u_1$ &  $A(-0.7060755\cos(2\pi x) - 0.0284767\sin(2\pi x))$ &  $A (-0.7201359 \cos(2\pi x) - 0.2482996 \sin(2\pi x))$ \\
 $u_2$ & $A (-0.7072239\cos(2\pi x) - 0.0002393\sin(2\pi x))$  &  $A (-0.4672883 \cos(2\pi x) - 0.3976914 \sin(2\pi x))$ \\     
 $u_3$ & $A (-0.7072145\cos(2\pi x) + 0.0025891\sin(2\pi x))$  &  $A (-0.1801365 \cos(2\pi x) - 0.3357016 \sin(2\pi x))$ \\

\end{tabular}
\end{minipage}
\end{table*} 

\subsubsection{DustyShock}

We publish the results of 9 runs of DustyShock Problem for 1, 2 and 3 dust fractions. In contrast to DustyWave models, in which each dust fraction has fixed stopping time, in DustyShock models each fraction has a fixed radius of particle $s_i$, which is used to define stopping time. For this purpose, we assumed that the dispersed phase interacts with gas in the Epstein mode \cite{Epstein1924}, material density of dust $\rho_{\rm s}=1$, and we obtained relation between the size of particle and local value of its stopping time: 
\begin{equation}
\label{eq:tstopi}
t_i = \frac{s_i \rho_{\rm s}}{c_{\rm s} \rho_{\rm g}}=\frac{s_i} {\displaystyle\sqrt{\frac{p}{\rho_{\rm g}}} \rho_{\rm g}}.    
\end{equation} 

In all calculations, we used the adiabatic exponent of gas $\gamma=1.4$. 

In models DS1-DS9, the problem is stated on interval [0,1]. Moreover, in all models we used the same initial data for gas, i.e., constant values of density, pressure, and internal energy, which undergo discontinuity in $x=0.5$:
$$[\rho_{{\rm g}}, p, e]_{L} = [1, 1, 2.5],$$
$$[\rho_{{\rm g}}, p, e]_{R} = [0.125, 0.1, 2].$$ 
Gas and dust fractions have zero velocity at the initial moment. The mass of SPH particle for each dust fraction is defined by gas density and    $\varepsilon_i$, which is specified identically both on the left and on the right of the discontinuity point. For each model, the value of each dust fraction $\varepsilon_i$ is presented in Table ~\ref{tab:runsetup}.   

To give the boundary condition (\ref{eq:boundaryDShock}), we placed motionless particles in domains $x \in [-2h,0]$ and $x \in [1,2h]$. In DS1-DS9, we used standard parameters of artificial viscosity $\alpha = 1$, $\beta = 2$, $\nu =0.1h$ in \cite{SPH}. 

\subsection{Polydisperse dust with short stopping time, asymptotic properties of the method}
\label{sec:smalltstop}

For the problems with small parameters, asymptotic preserving (AP) methods are efficient. According to the definition in \cite{Jin2010}, the AP method is the method, which reduces the error of numerical solution at constant time and space step as the small parameter of the problem decreases. To demonstrate this feature in practice, one needs to know the analytical solution of the problem for different stopping time values or to be able to evaluate the error of the solution by Runge’s technique. For DustyWave and DustyShock problems on wave propagation, we know reference solutions, but these solutions are approximate. For DustyWave, the reference solution is the solution of linearized problem and for DustyShock, the reference solution is the asymptotic solution of nonlinear problem independent of $t_i$. Moreover, for the Lagrangian methods, the applicability of error evaluation by Runge’s technique is not studied. Therefore, we will not study the property of the method according to definition. Instead, we will investigate if the method we proposed allows us to obtain solutions of acceptable precision without constraining time step and smoothing radius according to the small parameters of the problem. 

\begin{figure}[!t]
\centering
\includegraphics{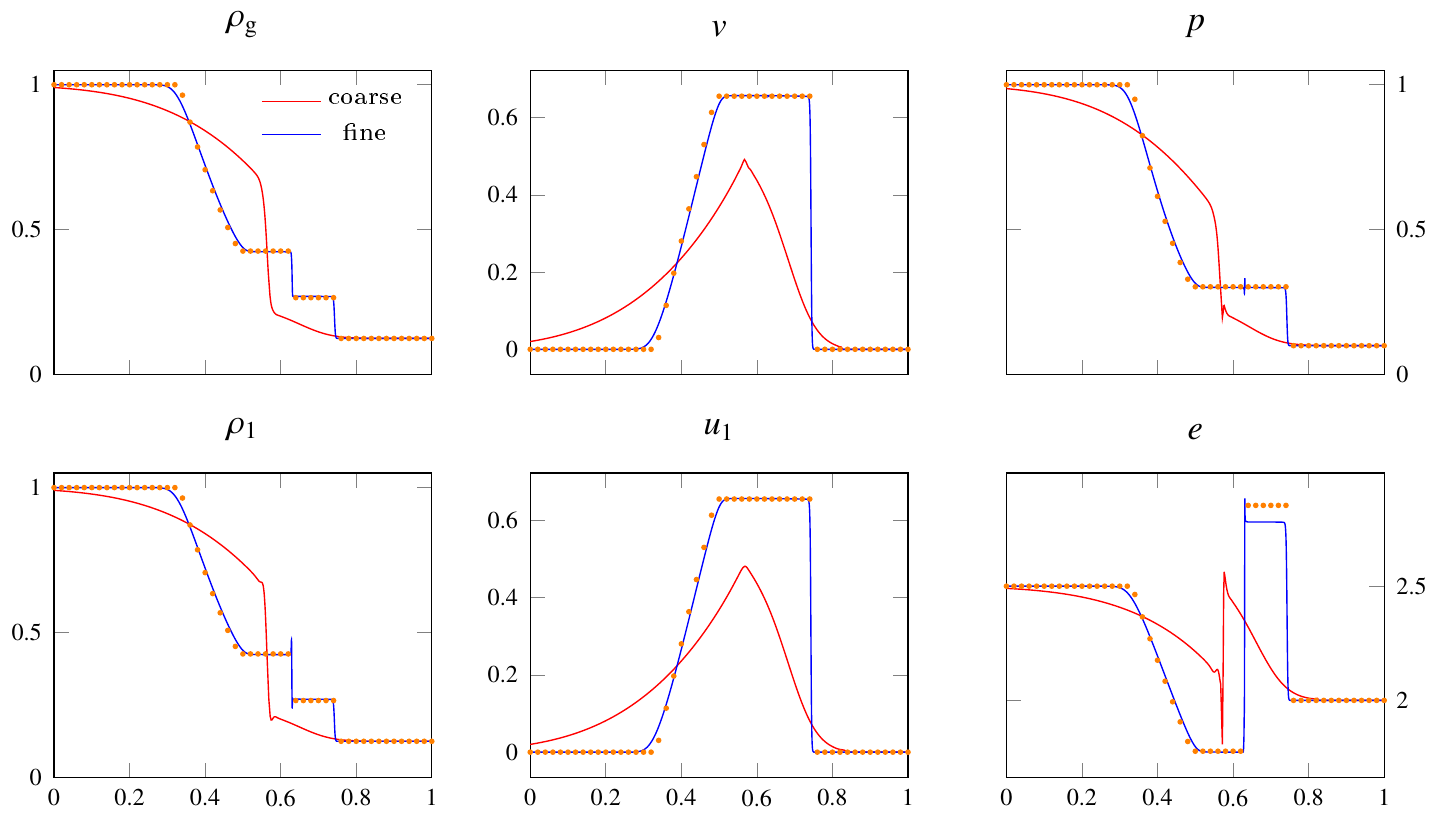}
\caption{The solution of shock-tube problem at time moment $T=0.2$ by the MK scheme (\ref{eq:MonKoch1995v})-(\ref{eq:MonKoch1995u}) for monodisperse dust. Dots are the reference solution and solid lines are numerical solutions. Red line is $N_{\rm SPH}=2100$, $h=0.01$, $\tau=5 \times 10^{-5}$ (model DS2). Blue line is  $N_{\rm SPH}=21000$, $h=0.001$, $\tau=5 \times 10^{-5}$ (model DS4). The MK scheme does not preserve asymptotics of the solution; therefore, for small stopping time values, to achieve an acceptable precision, it is necessary to take $\tau$ and $h$ basing on the value $t_{\rm stop}$.}
\label{fig:SmallGrainsMonoMK}
\end{figure}

\begin{figure}[!t]
\centering
\includegraphics{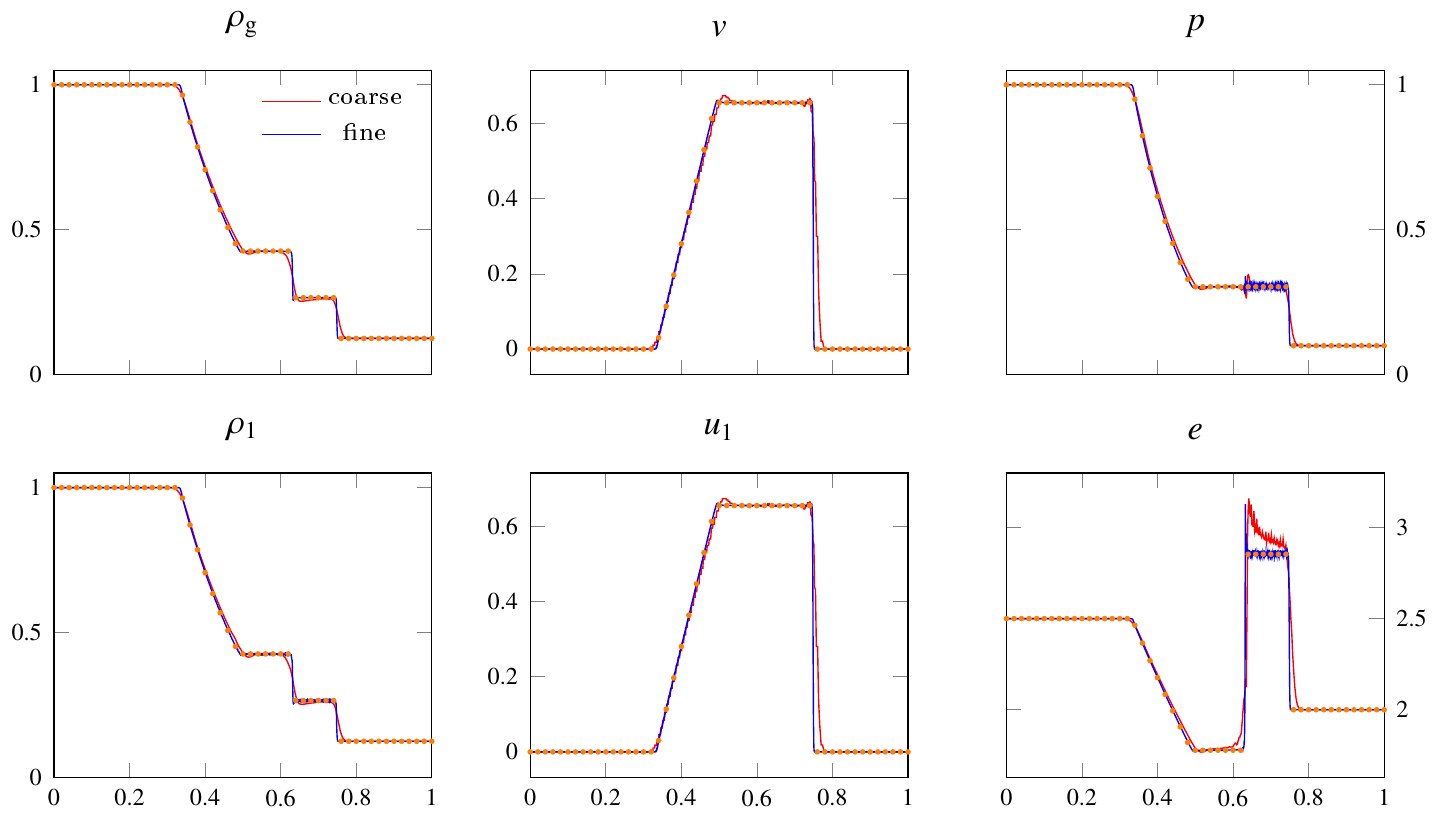}
\caption{The solution of shock-tube problem at time moment $T=0.2$ by the IDIC scheme with one dust fraction. Dots are the reference solution and solid lines are numerical solutions. Red line is $N_{\rm SPH}=2100$, $h=0.01$, $\tau=5 \times 10^{-3}$ (model DS1). Blue line is $N_{\rm SPH}=21000$, $h=0.001$, $\tau=5 \times 10^{-4}$ (model DS3). The IDIC scheme preserves asymptotics of the solution, therefore, for small stopping time values, an acceptable precision is achieved for $\tau$ and $h$ not depending on $t_{\rm stop}$, but satisfying the Courant condition for gas.}
\label{fig:SmallGrainsMonoIDIC}
\end{figure}

In order to make the term "acceptable precision" more definite, we use the methods of calculation of interphase interaction for monodisperse dust, which were systematically studied in \cite{LaibePriceAstroDrag,Booth2015}. Relating to the Monaghan-Kocharyan (MK) method \cite{MonaghanKocharyan1995}, we know empirical criterion (\ref{eq:h}), which associates smoothing radius and stopping time. If the smoothing radius does not satisfy (\ref{eq:h}), we observe a numerical dissipation of the solution caused by the calculation of intense interphase interaction \cite{LaibePrice2011Test}. This means that the MK method is not the AP method. We calculate DS2 of DustyShock problem with small stopping time by using the MK method with $h$, which does not satisfy (\ref{eq:h})\footnote{Our realization of the MK method is described in \ref{sec:MK}}. The level of dissipation can be evaluated visually in Fig.~\ref{fig:SmallGrainsMonoMK}, which shows numerical (red solid line) and reference (orange dots) solutions of this problem. It is seen that the spatial resolution of DS2 does not reproduce the wave structure properly. However, if we decrease the smoothing radius by a factor of 10 increasing the number of SPH particles (the calculation of DS4), the location of shock wave, contact discontinuity, and rarefaction wave will be shown correctly (Fig.~\ref{fig:SmallGrainsMonoMK}, blue line). Then we calculate this problem by using IDIC method and the same coarse (DS1) and fine (DS3) spatial resolution as in DS2 and DS4. It is seen from Fig.~\ref{fig:SmallGrainsMonoIDIC} that the calculations by using the IDIC method are much less sensitive to changes of $h$ and reproduce the wave structure correctly even at "coarse" smoothing radius. For polydisperse dust, we will consider the level of dissipation obtained by using the IDIC method with "coarse" resolution as "acceptable precision", and the level obtained by the MK method with "coarse" resolution will be considered as overdissipation. 

\begin{figure}[!t]
\centering
\includegraphics{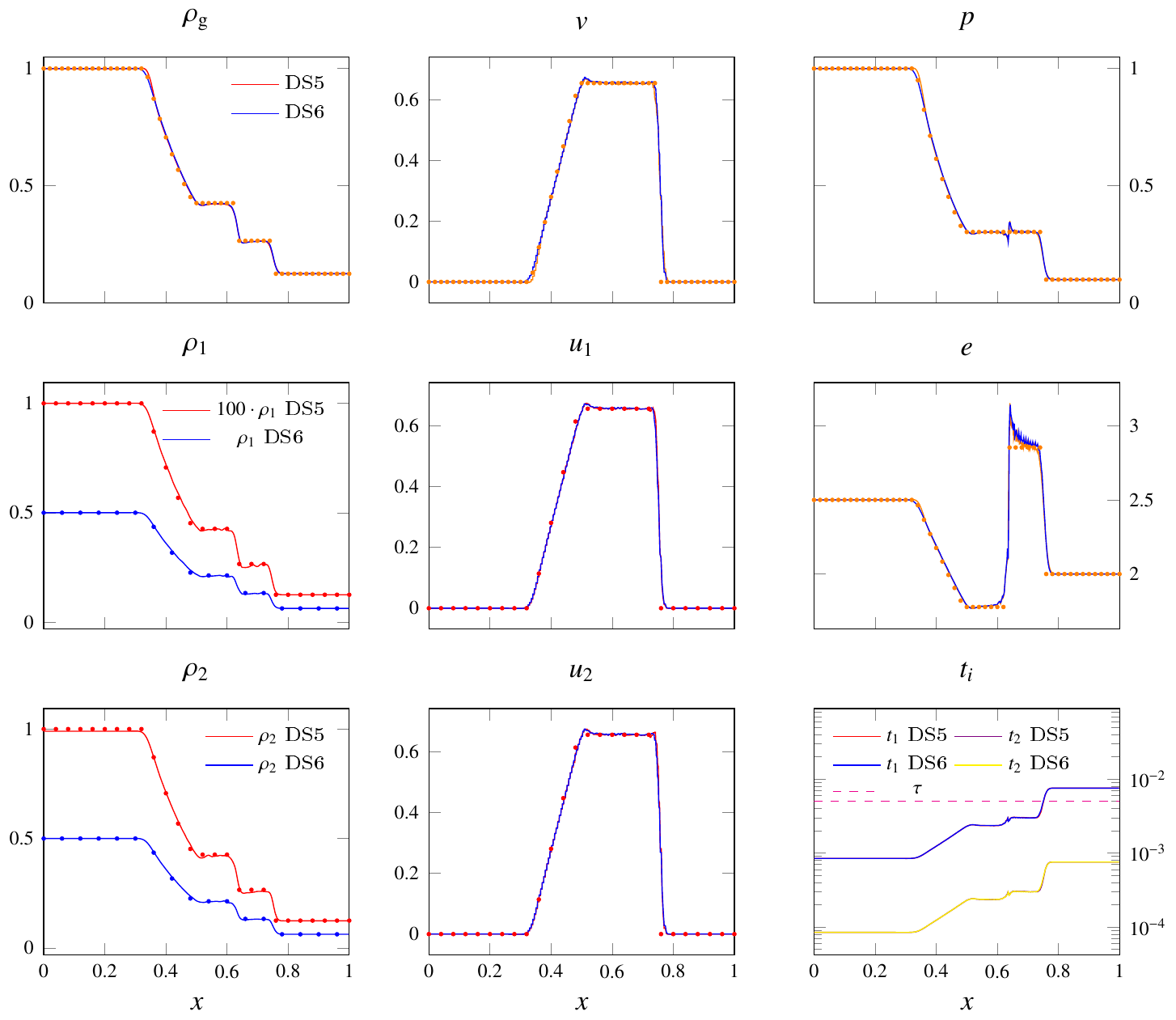}
\caption{The solution of shock-tube problem at time moment $T=0.2$ by the IDIC scheme for 2 fractions of polydisperse dust for small stopping time values. Dots are the reference solution and solid lines are numerical solutions. Red line is the model DS5 with $\varepsilon_1=0.01$, $\varepsilon_2=0.99$, blue line is the model DS6 with $\varepsilon_1=\varepsilon_2=0.5$. In both models $N_{\rm SPH}=1180$, $h=0.01$, $\tau=5 \times 10^{-3}$. The right-hand lower panel shows that the calculation with $\tau>\min(t_i)$.} \label{fig:9pan_SmallMultigrainWithTstop_2fraction}
\end{figure}

In models DS5 and DS6, we reduce the number of particles in comparison to "coarse" resolution of monodisperse models preserving all other calculation parameters of model DS1. We add the second dust fraction, which also has short stopping time, into the mixture. Let the dust to gas mass ratio of the finest dust be $\varepsilon_2=0.99$ in model DS5 and $\varepsilon_2=0.5$ in model DS6, but the total dust to gas mass ratio is the same in both models and equals 1. Due to (\ref{eq:csmulti}), we expect that $\rho_{\rm g}$, $v$, $u_i$, $p$, $e$ for these models will be the same, given $\varepsilon_i$ do not vary with time. The calculation results shown in Fig.~\ref{fig:9pan_SmallMultigrainWithTstop_2fraction} confirm our expectations: on the panels, which show the parameters of gas and dust velocity, the lines for DS5 and DS6 almost coincide. Moreover, comparing Fig.~\ref{fig:SmallGrainsMonoIDIC} and Fig.~\ref{fig:9pan_SmallMultigrainWithTstop_2fraction}, we see that passing from one to two dust fractions, the scale of numerical dissipation of solution remains acceptable. From the right-hand lower panel, it follows that stopping time of the finest dust in the model is within the range $10^{-4}-10^{-3}$. This means that the proposed method of calculation of interaction between phases allows us to find the solution even if the following conditions are not fulfilled: 
\begin{equation}
\label{eq:taumultigrain}
   \tau < \min(t_i),
\end{equation}
\begin{equation}
\label{eq:hmultigrain}
   h < \min(t_i) \cdot c_{\rm s}.
\end{equation}

\begin{figure}[!t]
\centering
\includegraphics{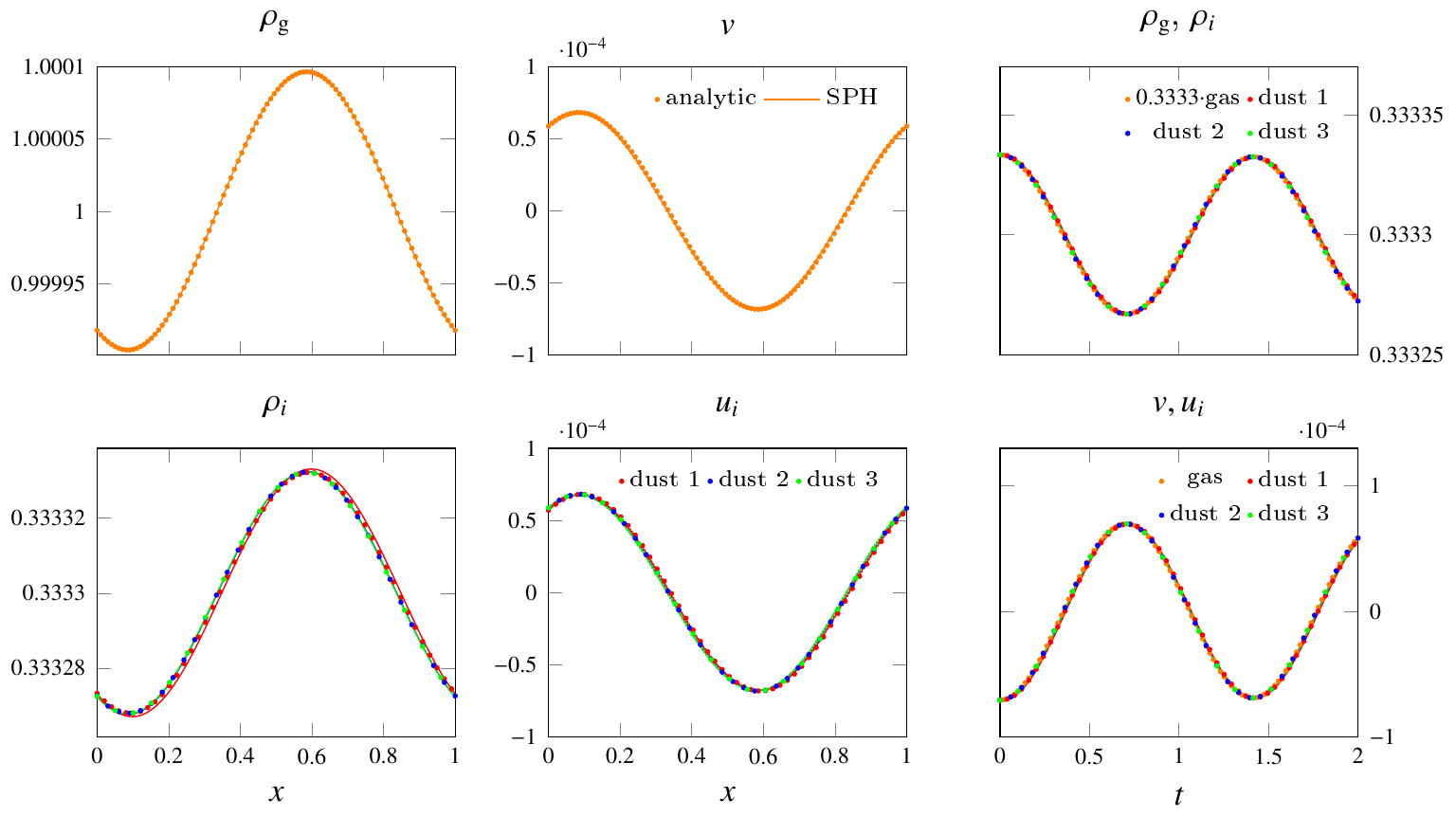}
	\caption{The solution of the problem on propagation of acoustic waves for small stopping time values ($t_1=10^{-2}$, $t_2=10^{-3}$, $t_3=10^{-4}$). Dots are the reference solution and solid lines are the calculation of model DW2 with $N_{\rm SPH}=600$, $h=0.01$, $\tau=5 \times 10^{-3}$. Gas parameters are shown in orange, three dust fractions are shown in red, blue, and green. Four left-hand panels are the dependence of the solution on the coordinate at time moment $T=2$. Two right-hand panels demonstrate the dependence of the solution on time at point $x=0$. The damping of waves is absent.}
\label{fig:Wave_SmallMultigrain_3fraction}	
\end{figure}

We confirm this conclusion taking into consideration the DustyWave problem with 3 fractions of fine dust. In comparison to DS5 and DS6, the DW2 model has additional coarser dust fraction, which still moves at almost gas velocity.  The mass fraction $\varepsilon_i$ of all three fractions in DW2 is identical. Due to the fact that the solution of DW2 is continuous, we decrease the number of SPH particles to 600, but other numerical parameters remain unchanged. The results of calculation of model DW2 are shown in Fig.\ref{fig:Wave_SmallMultigrain_3fraction}. The left-hand and central panels in the figure demonstrate density and velocity of gas and dust at time moment $T=2$. It is seen that all the dust fractions have velocity values close to each other and to gas velocity. The right-hand panels in the figure show that initial amplitudes of perturbations do not vary with time. These calculations confirm the conclusion on non-damping of sound waves for small $t_i$ known from the dispersion relation analysis. Moreover, the coincidence of numerical and reference solutions shows that the method has an acceptable level of numerical dissipation even under violation of conditions (\ref{eq:hmultigrain})-(\ref{eq:taumultigrain}).  

%\textcolor{magenta}{Чтобы предварительно оценить применимость метода в пространственно трехмерном случае, мы рассчитали модель DW3 с существенно более грубым пространственным разрешением. Мы сохранили все физичекие параметры модели DW2, но вместо 600 частиц на длину волны в DW2 мы использовали 30 частиц на длину волны в DW3. Это означает, что в DW3 необходимо создать профиль возмущения малой амплитуды $A=10^{-4}\rho_0$, используя модельные частицы большой массы $1/30 \int_0^1 \rho_{\rm g}(x) dx$. Поэтому профиль волны в таком расчете будет чувствителен к положению модельных частиц. Эта особенность задачи налагает свои ограничения на радиус сглаживания и шаг по времени. Поэтому в DW3 мы не смогли пропорционально увеличить в 20 раз все расчетные параметры, а использовали $h=0.1$ и $\tau=0.005$. При таких параметрах в расчете трения средние скорости газа и пыли в ячейках вычисляются по 1-2 частицам, поэтому при перемещении частицы в соседнюю ячейку следует ожидать значительных осцилляций в скоростях. The results of DW3 simulation are shown in Fig.~\ref{fig:Wave_SmallMultigrain_3fraction_lowResolution}. Видно, что использование экстремально грубого разрешения не приводит к уменьшению амплитуды волны, однако переводит исходно монотонный профиль скорости в осцилирующий. Это означает, что расчет движения волн на длительное время может потребовать коррекции скорости искусственной вязкостью. }

\begin{figure}[!t]
\centering
\includegraphics{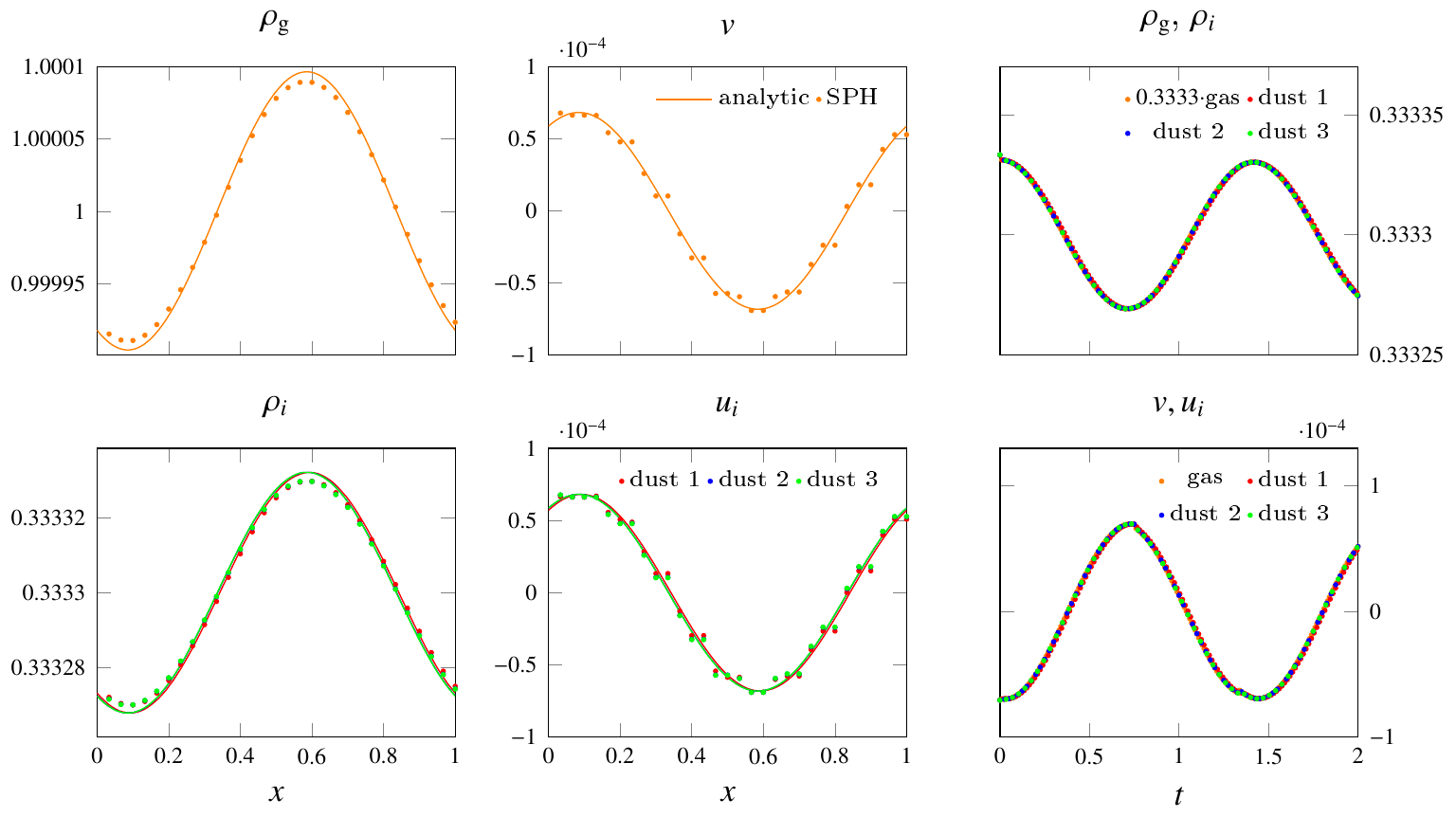}
	\caption{The solution of the problem on propagation of acoustic waves for small stopping time values ($t_1=10^{-2}$, $t_2=10^{-3}$, $t_3=10^{-4}$). Solid lines are the reference solution and dots are the results of calculation of low-resolution model DW3 with $N_{\rm SPH}=30$, $h=0.1$, $\tau=5 \times 10^{-3}$. Gas parameters are shown in orange, three dust fractions are shown in red, blue, and green. Four left-hand panels are the dependence of the solution on the coordinate at $T=2$. Two right-hand panels demonstrate the dependence of the solution on time at $x=0$. The initial monotonic velocity profile becomes oscillating.}
\label{fig:Wave_SmallMultigrain_3fraction_lowResolution}	
\end{figure}

To pre-evaluate the applicability of the method to three-dimensional simulations, we compute the DW3 model with a significantly coarser spatial resolution. We keep all the physical parameters of the DW2 model, but instead of 600 particles per wavelength in DW2, we use 30 particles per wavelength in DW3. This means that in DW3 it is necessary to create a small amplitude wave profile $A=10^{- 4} \rho_0$ using model particles of large mass $ 1/30 \int_0^1 \rho_{\rm g}(x)dx$. Therefore, the wave profile in DW3 is much more sensitive to the position of the model particles than in DW2. This feature of the problem imposes its own constrains on the smoothing length and time step. Therefore, in DW3, we can't just scale all the numerical parameters by 20 times, but instead we use $h=0.1$ and $\tau = 0.005$. To compute drag-in-cell terms with such parameters, the mean velocities of gas and dust in the cells are averaged over 1–2 particles; therefore, when a particle moves to the next cell, significant oscillations in the velocities is expected. The results of DW3 simulation are shown in Fig.~\ref{fig:Wave_SmallMultigrain_3fraction_lowResolution}. One can see that the use of an extremely coarse resolution does not lead to wave overdamping, but transforms the initially monotonic velocity profile into an oscillating one. This means that the calculation of wave propagation for a long time may require correction of the oscillating velocity with artificial viscosity.

All the results point to the fact that the method of calculation of interphase interaction for polydisperse dust "inherited" asymptotic features of the method for monodisperse dust. 

\subsection{Polydisperse dust with long and multi-scale stopping time values}
\label{sec:multiscaletstop}

\begin{figure}[!t]
\centering
\includegraphics{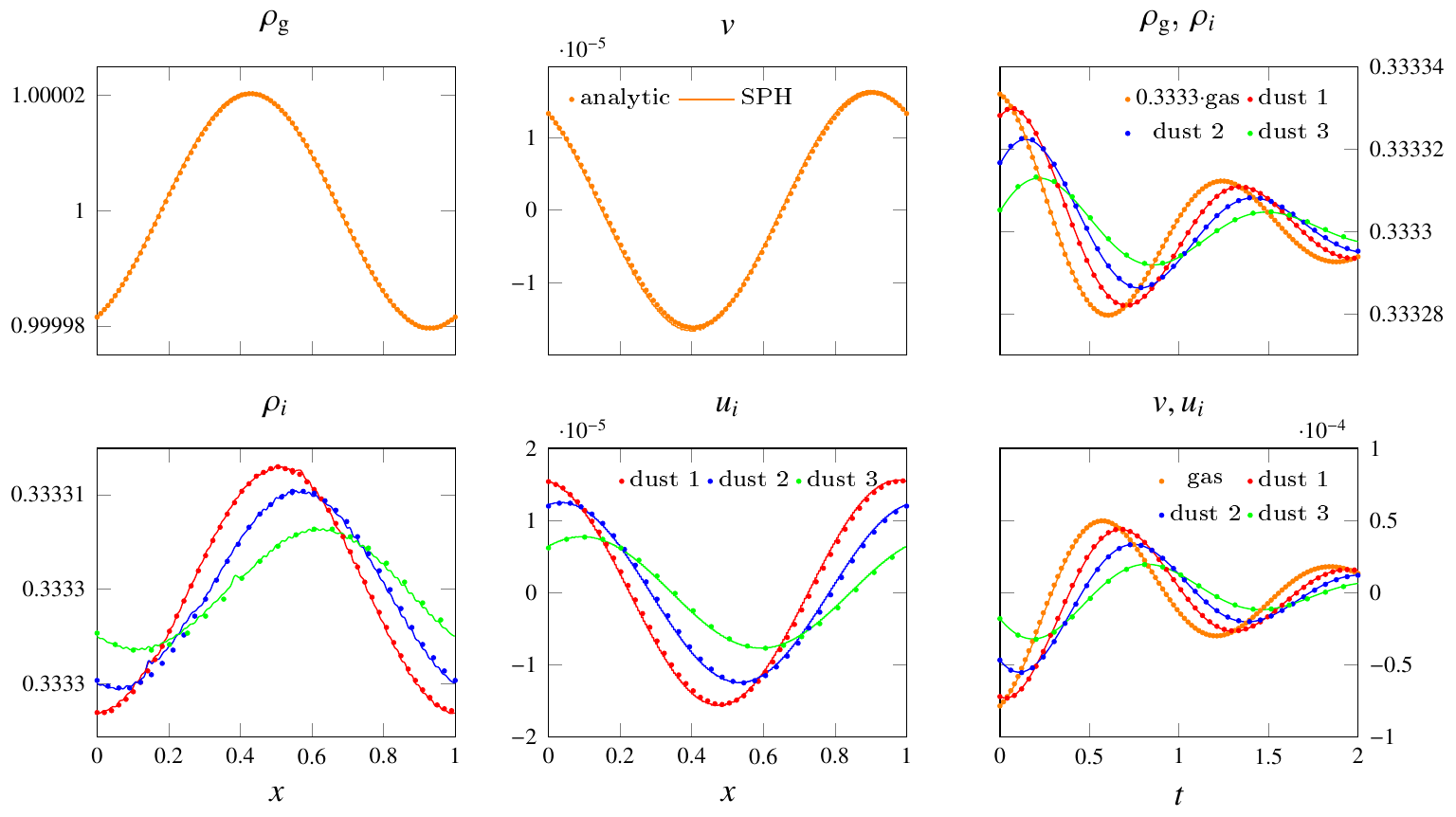}
\caption{The solution of the problem on propagation of acoustic waves for 3 dust fractions with large stopping time values ($t_1=0.1$, $t_2=0.2$, $t_3=0.4$). Dots are the reference solution, solid lines are the calculation of model DW1 with $N_{\rm SPH}=600$, $h=0.01$, $\tau=5 \times 10^{-3}$. Gas parameters are shown in orange, three dust fractions are shown in red, blue, and green. Four left-hand panels are the dependence of the solution on the coordinate at time moment $T=2$. Two right-hand panels demonstrate the dependence of the solution on time at point $x=0$. The physical damping of waves takes place.}
 
\label{fig:Wave_BigMultigrain_3fraction}
\end{figure}

Before we pass to the most challenging problem related to reproducing shock solutions in the gas-polydisperse dust medium with multi-scale stopping time values, we make sure that the proposed method is successfully used to reproduce large stopping time values. For this purpose, we calculate the propagation of acoustic waves by the model DW1 with three fractions of coarse dust. The numerical parameters of DW1 coincide with those of DW2. The results of calculation are presented in Fig.~\ref{fig:Wave_BigMultigrain_3fraction} where the numerical solution is shown by solid lines and the reference solution is demonstrated by dots. It is seen that for all variables the numerical and reference solutions coincide. In the central panels of the figure, it is seen that in contrast to the model with short stopping time values, the velocity values of gas and dust fractions differ. In the right-hand panels of Fig.~\ref{fig:Wave_BigMultigrain_3fraction}, where the dependence of the solution on time at point $x=0$ is shown, it is seen that in contrast to the model DW2, in DW1 the damping of perturbation amplitude with time takes place. Such a behavior of the solution agrees with the results of dispersion relation analysis for mono- and polydisperse dust \cite{LaibePrice2011Test,BenitezLlambay2019}. 
We also note that in the model DW1 for large stopping time values the density values of dust fractions are less smooth than in the model DW2. A higher noise level in DW1 is resulted from that the velocity of coarse dust fractions is significantly different from gas velocity while those for fine dust is almost indistinguishable. However, due to the fact that in DW1, the averaging of the velocity values of gas and dust fractions is made with respect to only 3 model particles, the obtained noise level can be considered as acceptable. 

In conclusion, we study the ability of the method to reproduce the dynamics of polydisperse particles for small and large stopping time values under violation of conditions (\ref{eq:taumultigrain})-(\ref{eq:hmultigrain}) in DustyShock problem. In contrast to DustyWave problem, the reference solution of DustyShock for the mixture of coarse and fine particles is absent. We will judge about the quality of the obtained numerical solution by using calculations with diminished $h$ and $\tau$.  

\begin{figure}[!t]
\centering
\includegraphics{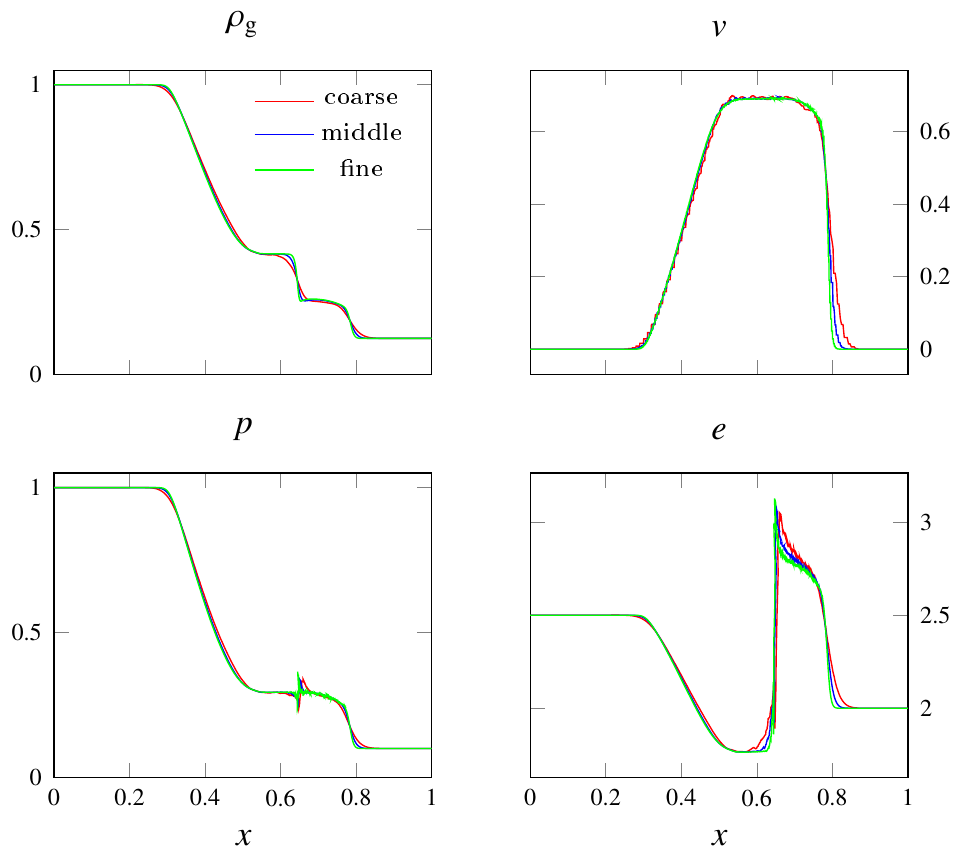}
\caption{ The solution of shock-tube problem at time moment $T=0.2$ for polydisperse dust with multi-scale stopping time values. The calculations with different numerical resolution are illustrated by different colors. The model DS7 with $N_{\rm SPH}=1180$, $h=0.02$, $\tau=5 \times 10^{-3 }$ is red, the model DS8 with $N_{\rm SPH}=2360$, $h=0.01$, $\tau=2.5 \times 10^{-3}$ is blue, the model DS9 with $N_{\rm SPH}=7086$, $h=0.005$, $\tau=1.25 \times 10^{-3}$ is green. The panels show gas parameters.}
\label{fig:3fractionDustyShock_gas}
\end{figure}

\begin{figure}[!t]
\centering
\includegraphics{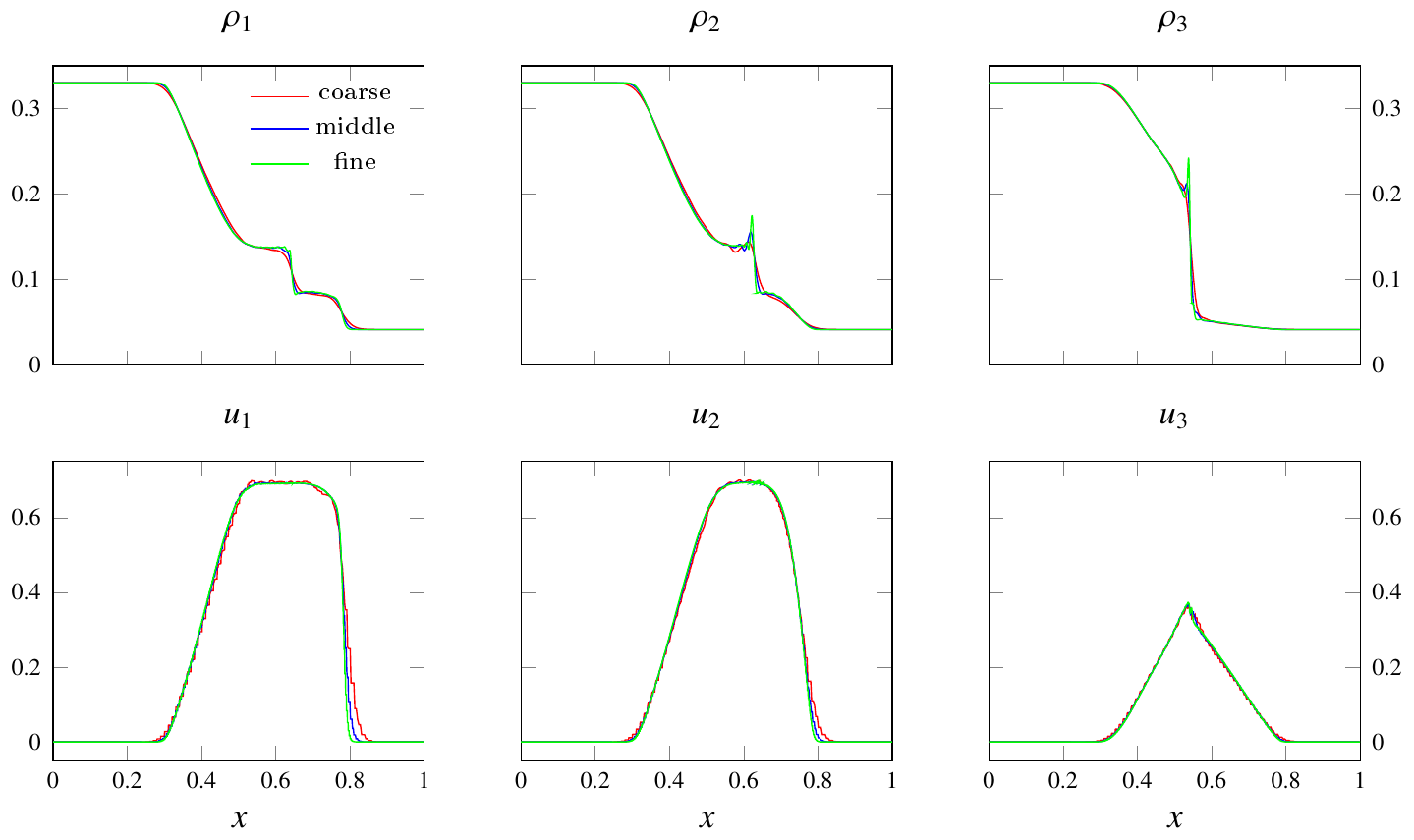}
\caption{The solution of shock-tube problem at time moment $T=0.2$ for polydisperse dust with multi-scale stopping time values. The calculations with different numerical resolution are illustrated by different colors. The model DS7 with $N_{\rm SPH}=1180$, $h=0.02$, $\tau=5 \times 10^{-3}$ is red, the model DS8 with $N_{\rm SPH}=2360$, $h=0.01$, $\tau=2.5 \times 10^{-3}$ is blue, the model DS9 with $N_{\rm SPH}=7086$, $h=0.005$, $\tau=1.25 \times 10^{-3}$ is green. The upper panels show the density of dust fractions, the lower ones show the velocity of dust fraction.}
\label{fig:3fractionDustyShock_dust}
\end{figure}

We consider the model DS7 with three dust fractions. Let a finely dispersed fraction be of size $s_1=10^{-3}$, and the second and the third fractions be coarser than the first one by a factor of 10 and 100, respectively. Assume that the dust to gas mass ratio $\varepsilon_i$ of each fraction is identical. For base model DS7, we take $\tau$ and $h$, for which conditions (\ref{eq:taumultigrain})-(\ref{eq:hmultigrain}) are violated. In the models DS8 and DS9, we will increase the number of particles and diminish the smoothing radius. The results of calculations of models DS7-DS9 are shown in Fig.~\ref{fig:3fractionDustyShock_gas}
and Fig.~\ref{fig:3fractionDustyShock_dust}. It is seen from Fig.~\ref{fig:3fractionDustyShock_gas} that the positions of shock wave fronts, contact discontinuity, and rarefaction wave in gas are reproduced similarly with different $\tau$ and $h$. Moreover, the velocity of all three dust fractions and the position of wave fronts are reproduced similarly in all models. Features related to resolution arise in the density of medium and coarse dust in the vicinity of contact discontinuity for these fractions where the density bump and the observable wave structure take place. From central and right-hand panels in Fig.~\ref{fig:3fractionDustyShock_dust}, it follows that the higher the spatial resolution is, the more localized and expressed this bump is. The nature of this bump and waves requires individual study, which is beyond the scope of this paper. However, the behavior of this bump observed on a qualitative level corresponds to the features of dispersion relation (\ref{1.12}) for gas and three dust fractions. It follows from this relation that short-wave perturbations undergo strong dispersion. 

In general, the calculations of polydisperse dust dynamics for large and multi-scale stopping time values show that the proposed method of drag calculation does not introduce dissipation to the solution and the density value of dust fractions are more sensitive to appearance of perturbations of both physical and numerical nature.

\section{Conclusions}
\label{sec:resume}

In the paper, we present the IDIC (Implicit-Drag-in-Cell) method of calculation of interphase interaction to model the flows of two-phase polydisperse media based on smoothed particles hydrodynamics. In polydisperse media, stopping time values of dispersed inclusions differ by the orders of magnitude and can be small parameters of the problem. The numerical solution of the problems with a small parameter by the SPH method can require that the smoothing radius and time step be defined by this parameter. This means that the closer the value of small parameter to zero is, the more computational costs on the problem solution are. The IDIC method of calculation of interphase interaction proposed in the paper solves this problem, i.e., it allows us to model the dynamics of two-phase polydisperse media with the same numerical resolution as the dynamics of pure fluid.

The IDIC method has the Lagrange-Euler nature, i.e., here interphase interaction is calculated by using Euler's grid and other forces are calculated by using the Lagrangian approach. The use of Euler's grid enables us to provide such an implicit approximation of terms responsible for momentum exchange between the carrier phase and dispersed inclusions that its realization does not require iteration process and a local law of conservation of momentum be fulfilled with machine precision. In the IDIC method, to find the velocity of carrier and dispersed phases at the next time step, we need to solve the system of linear equations with matrix of special type. Owing to a special structure of the matrix, the solution can be found by a direct method in $O(N^2)$ operations instead of standard $O(N^3)$, where $N$ is the number of dust fractions. 

We showed that the IDIC method does not introduce an additional dissipation into the numerical solution, which is related to the calculation of intense interphase interaction. On the contrary, it introduces dispersion whose intensity depends on problems, stopping time values, and resolution. We found that in case of long stopping time, it will be clearly seen in density of dust fraction. But in case of short stopping time and low resolution, the dispersion will manifest itself in oscillation of velocity profiles.

Overall, the combination of standard SPH (particle-particle $PP$) way of computing hydrodynamic forces and cell averaging (particle-mesh $PM$) way of computing drag results in counteraction of multi-directional numerical effects: dissipation and dispersion. We found that such $PPPM$ combination appears beneficial for the problems with several stiff relaxation terms since it provides computationally cheap and synergistic solution both to dissipation (also known as overdamping) and timestepping problem. Such solution requires extra memory and loses pure Lagrangian nature of the method. However, there is no reason to assume, that IDIC will automatically provide benefits in struggle with other computational problems of multifluid smoothed particle hydrodynamics. In particular the ability of IDIC to deal with dust particle clumping in regions where gas is under-resolved will be studied in future. 

%We also demonstrated that in the calculation of dynamics of gas-dust media for large stopping time values in the solution obtained by this method, there is a dispersion, which is more pronounced in the density of dust fractions.  

An attractive feature of the IDIC method is that it does not specify requirements on grid structure, on which drag is calculated. So the grid cells can be of arbitrary shape. However, in the presented version of the method, it is impossible to compute drag in half empty cells. It means that at least one model particle of each dust fraction must fall into every cell where gas particle is. In the nearest future, we plan to eliminate this restriction and to extend the method to three-dimensional problems. 

\section{Acknowledgment}

This study was funded by the Russian Science Foundation grant number 19-71-10026. We thank anonymous reviewers for their constructive suggestions and helpful comments, we are grateful to Olga Drozhzhina for her assistance with language editing.

\appendix 

\section{The calculation of interphase interaction via drag by the Monaghan-Kocharyan method}
\label{sec:MK}

The solutions illustrated in Fig.~\ref{fig:SmallGrainsMonoMK} were obtained by the explicit Monaghan-Kocharyan scheme \cite{MonaghanKocharyan1995} for monodisperse dust:

\begin{equation}
\displaystyle\frac{v^{n+1}_a-v^n_a}{\tau}= - \sum_b m_{\rm g} \left(\frac{P_b}{(\rho^n_{b,\rm{g}})^2} + \frac{P_a}{(\rho^n_{a,\rm g})^2} + \Pi_{ab} \right) \bigtriangledown_a W^{n}_{ab} 
-m_{\rm d} \sum_l \frac{K_{al}}{\rho^n_{a, \rm g} \rho^n_{l,\rm d}} \frac{(v_a^{n} - u_l^{n},r_{la})}{r_{la}^2+\eta^2}r_{la}W^{n}_{la},
\label{eq:MonKoch1995v}
\end{equation}
\begin{equation}
\displaystyle\frac{u^{n+1}_l-u^n_l}{\tau}= m_{\rm g} \sum_a \frac{K_{al}}{\rho^n_{a, \rm g} \rho^n_{l,\rm d}} \frac{(v_a^{n} - u_l^{n},r_{la})}{r_{la}^2+\eta^2}r_{la}W^{n}_{la},
\label{eq:MonKoch1995u}
\end{equation}

\begin{equation}
K_{al}=\displaystyle\frac{\rho^n_{l,\rm d} \rho^n_{a, \rm g} c^n_{a,\rm s}}{s_l^n \rho_{\rm s}},
\end{equation}

where $r_{la} = r_l - r_a$, $\eta$ is a clipping constant, $\eta^2 = 0.001 h^2$, $\rho^n_{l,\rm d}$ is dust density in $l$-numbered particle , $m_{\rm d}$ is a mass of model dust particle, $K_{al}$ is a drag coefficient between $a$-numbered gas particle and $l$-numbered dust particle. In the calculations, we used constant particle size $s_l^n=10^{-4}$ and constant material density of dust particles $\rho_{\rm s}=1$.

\section{The type of gas-polydisperse dust equations in isothermal case}
\label{sec:type}

The theory of existence and uniqueness for smooth and discontinuous solutions of initial-boundary value problems is widely developed for hyperbolic systems of partial differential equations. It is known that the equations of mechanics of two-phase gas-monodisperse dust medium can be of hyperbolic and composite type (see, for example, \cite{FedorovBedarev2018}). The composite type of the system is expected to result in additional complications in its numerical solution. We make sure that the solved problem of constructing AP method is unrelated to the properties of equations and mathematical statement of the problem, but it has only a numerical nature. We show that one-dimensional system of isotermal gas-polydisperse dust equations is hyperbolic. 

We rewrite the system (\ref{eq:contSPHgas})-(\ref{eq:motionSPHdust}) and (\ref{eq:pisothermal}) as follows
\begin{equation}
\label{eq:systemDiv}
\displaystyle\frac{\partial \Phi}{\partial t}+M\frac{\partial \Phi}{\partial x}=\Psi,
\end{equation}

\begin{equation}
   \Phi=\begin{pmatrix}
   \rho_{\rm g}\\
   v\\
   \rho_1\\
   u_1\\
   ...\\
   \rho_N \\
   u_N
   \end{pmatrix},
\quad \quad   
   \Psi=\begin{pmatrix}
   0\\
   f_{\rm g}\\
   0\\
   f_1\\
   ...\\
   0\\
   f_N
   \end{pmatrix},
\quad \quad   
   M=\begin{pmatrix}
   v & \rho_{\rm g} & 0 & 0 &...& 0 & 0 \\
   \displaystyle\frac{c_{\rm s}}{\rho_{\rm g}} & v & 0 & 0 &...& 0 & 0\\
   0 & 0 & u_1 & \rho_1 &...& 0 & 0\\
   0 & 0 & 0 & u_1 & ...& 0 & 0\\
   ...&...&...&...&...&...&...&\\
   0 & 0 & 0 & 0 &...& u_N & \rho_N \\
   0 & 0 & 0 & 0 & ...& 0 & u_N 
   \end{pmatrix}.    
\end{equation}
Then it is easily seen that characteristic equation of matrix $M$ 
\begin{equation}
(u_N-\lambda)^2 \cdot...\cdot(u_1-\lambda)^2\left((v-\lambda)^2 -c_{\rm s}^2 \right)=0    
\end{equation}
has only real roots
\begin{equation}
\lambda_{2j-1,2j}=u_j, \quad j=1,..N, \quad \lambda_{2N+1,2N+2}=v \pm c_{\rm s}.  
\end{equation}

\section{Derivation of reference solution for isothermal gas-polydisperse dust system}
\label{sec:dustyWaveMultigrain}

The reference solution of DustyWave problem presented in Fig.~\ref{fig:Wave_SmallMultigrain_3fraction} and Fig.~\ref{fig:Wave_BigMultigrain_3fraction} was obtained by using our program DustyWaveMulti.sci. At the input, this code takes the parameters of gas-dust medium and size of domain $x \in [0,L]$, $t \in [0,T]$, in which we seek the solution. At the output, we obtain the analytical representation for $\rho_{\rm g}(x)$, $\rho_i(x)$, $v(x)$, $u_i(x)$ at time moment $t=0$, and numerical values $\rho_{\rm g}(x)$, $\rho_i(x)$, $v(x)$, $u_i(x)$ at time moment $t=T$, and numerical values $\rho_{\rm g}(t)$, $\rho_i(t)$, $v(t)$, $u_i(t)$ at point $x=x_0$.  

The method of obtaining the solution realized in DustyWaveMulti.sci is presented in \cite{BenitezLlambay2019}. In what follows, we describe the main steps of this method. The system (\ref{eq:DWS1})-(\ref{eq:DWS4}), which describes the dynamics of isothermal gas with $N$ disperse dust fractions, is nonlinear. The stationary solution (\ref{eq:steadySolDW}) satisfies this system. We linearize the system (\ref{eq:DWS1})-(\ref{eq:DWS4}) in the vicinity of (\ref{eq:steadySolDW}) and seek the solution of linear system. Under small deviations from stationarity, the solution of initial nonlinear and linearized systems will be close. Further, we consider the main steps for finding such a linearized solution.

Assume that the solution (\ref{eq:DWS1})-(\ref{eq:DWS4}) takes the form $\rho_{\rm g}(x,t)=\rho_{\rm g}^0+\delta\rho_{\rm g}(x,t)$, $\rho_j(x,t)=\rho^0_j+\delta\rho_j(x,t)$, and $v(x,t)=\delta v(x,t)$, $u_j(x,t)=\delta u_j(x,t)$. 
Substituting it into the system (\ref{eq:DWS1})-(\ref{eq:DWS4}) 
and neglecting the terms of the second order of smallness, we obtain the linear system with $2N+2$ differential equations 
\begin{align}\label{1.4}
   \frac{\partial\delta\rho_{\rm g}}{\partial t}+\rho_{\rm g}^0 \frac{\partial\delta v}{\partial x} &= 0,\\ \label{1.5}
   \frac{\partial\delta\rho_j}{\partial t}+\rho^0_j \frac{\partial\delta u_j}{\partial x} &= 0,\\ \label{1.6}
   \frac{\partial\delta\nu_g}{\partial t}+ \sum_{j=1}^{N}\frac{\varepsilon^0_j}{t_j} (\delta v-\delta u_j)+\frac{c^2_s}{\delta\rho_{\rm g}^0} \frac{\partial\delta\rho_g}{\partial x}&= 0,\\ \label{1.7}
   \frac{\partial\delta u_j}{\partial t}+\frac{1}{t_j}(\delta u_j-\delta v)&= 0,
   \end{align}
where the subscripts of dust particles $j=1,\cdots,N$ and $\varepsilon^0_j=\rho^0_j/\rho_{\rm g}^0$.

To solve this system, represent each of sought functions (\ref{1.4})-(\ref{1.7}) as 
\begin{equation}
\label{eq:perturbation}
\delta f=\delta\hat{f}e^{ikx-\omega t},
\end{equation}
where $k$ is a real wave number, $\omega$ is a complex frequency. The substitution of (\ref{eq:perturbation}) into \eqref{1.4} - \eqref{1.7} yields the system of linear equations for perturbations $\delta\hat{f}$ 
\begin{align}\label{1.8}
- \omega\delta \hat{v}+ \sum_{j=1}^{N}\frac{\varepsilon^0_j}{t_j} (\delta \hat{v}-\delta \hat{u}_j)+ik\frac{c^2_s}{\rho_{\rm g}^0}\delta\hat{\rho}_{\rm g} &= 0,\\ \label{1.9}
- \omega\delta \hat{u}_j+\frac{1}{t_{j}}(\delta \hat{u}_j-\delta \hat{v})&= 0,\\ \label{1.10}
-\omega\delta\hat{\rho}_{\rm g}+ik \rho_{\rm g}^0\delta \hat{v}&= 0,\\ \label{1.11}
-\omega\delta\hat{\rho}_j+ik \rho^0_j\delta \hat{u}_j&= 0. 
\end{align}
The system \eqref{1.8}-\eqref{1.11} has nontrivial solutions if its determinant is equal to 0. Setting to zero the determinant of the system 

   \begin{equation}\begin{vmatrix}
	-\omega+\sum_{j=1}^{N}\frac{\varepsilon^0_j}{t_{j}}& -\frac{\varepsilon^0_1}{t_{1}}& \ldots &-\frac{\varepsilon^0_{N}}{t_{N}}& \frac{c^2_sik}{\rho_{\rm g}^0}& 0& \ldots& 0\\
	-\frac{1}{t_{1}}& -\omega+\frac{1}{t_{1}}& \ldots& 0&  \ldots& \ldots& \ldots& 0\\
	\vdots& \vdots& \vdots& \vdots& \vdots& \vdots& \vdots& \vdots \\
	-\frac{1}{t_{N}}& 0& \ldots& -\omega+\frac{1}{t_{N}}& 0& \ldots&\ldots& 0\\
	ik&0& \ldots& 0& -\frac{\omega}{\rho_{\rm g}^0}& 0& \ldots& 0\\
	0& ik& 0& \ldots& 0& -\frac{\omega}{\rho^0_1}& \ldots& 0\\
	\vdots& \vdots& \vdots& \vdots& \vdots& \vdots& \ddots&\vdots \\
	0& \ldots&\ldots& ik& 0& \ldots&\ldots& -\frac{\omega}{\rho^0_{N}}\\
   \end{vmatrix}\end{equation}

we obtain the dispersion relation, which relates frequency and wave number 
\begin{equation}\label{1.20}
P(\omega) \equiv \displaystyle\frac{(-1)^N \omega^N}{\prod_{j=1}^N \rho^0_j}\left(\omega^2(\prod_{j=1}^{N}(1-\omega t_{j})+\sum_{j=1}^{N}\varepsilon^0_j\prod_{p\ne j}(1-\omega t_{p}))+\omega^2_{\rm s}\prod_{j=1}^{N}(1-\omega t_j)\right)=0,
\end{equation} 
where $\omega_{\rm s}=kc_{\rm s}$. 
By virtue of the fact that $P(\omega)$ is the polynomial of degree $(2N+2)$, Eq.~(\ref{1.20}) has $(2N+2)$ roots. The root $\omega=0$ has multiplicity $N$. Due to (\ref{eq:perturbation}), zero perturbation corresponds to root $\omega=0$, damping waves correspond to real and positive values of $\omega$, waves with increasing amplitude correspond to real negative values, and the periodic solutions correspond to purely imaginary values. We study the remained $(N+2)$ nonzero roots of the polynomial $P(\omega)$, which coincide with the roots of polynomial 
\begin{equation}
  \Tilde{P}(\omega)\equiv \omega^2(\prod_{j=1}^{N}(1-\omega t_{j})+\sum_{j=1}^{N}\varepsilon^0_j\prod_{p\ne j}(1-\omega t_{p}))+\omega^2_{\rm s}\prod_{j=1}^{N}(1-\omega t_j).   
\end{equation}
 For this purpose, we rewrite (\ref{1.20}) in the form convenient for analysis: 
\begin{equation}\label{1.12}
F(\omega,\omega_{\rm s})\equiv\omega^2(1+\sum_{j=1}^{N}\frac{\varepsilon^0_j}{1-\omega t_{j}})+\omega^2_{\rm s}=0.
\end{equation}

We first note that $F$ is positive for $\omega <\min(t^{-1}_j)$. Moreover, one-sided limits $F$ hold the condition
 \begin{equation}\label{1.21}
\lim_{\omega \to  t^{-1}_j \pm} F(\omega,\omega_{\rm s})=\mp \infty,
\end{equation}
i.e., the function changes sign when passing through a singular point $\omega=t_j^{-1}$. This implies that on the left from $\min(t^{-1}_j)$, the polynomial $\Tilde{P}(\omega)$ has no real roots, but between two neighboring points, there is at least one real positive root. The total number of such roots is $(N-1)$. Moreover, since
\begin{equation}\label{1.22}
\lim_{\omega \to \infty} F(\omega,\omega_{\rm s})=+\infty,
\end{equation}
there is one more real positive root of the polynomial $\Tilde{P}(\omega)$ located on the right from the last singular point of $F$. This means that polynomials $P(\omega)$ and $\Tilde{P}(\omega)$ have $N$ different real positive roots. We find all these roots by a numerical dichotomy method and denote $\omega_j, \quad j=1,..N$.
 
The remained two roots of the polynomial $\Tilde{P}(\omega)$ $\omega_{N+1}$ and $\omega_{N+2}$ can be expressed through found roots by using Vieta's formula 
\begin{align}
 \omega_{N+1} \omega_{N+2} &= (-1)^{N+2}a_0\prod_{j=1}^{N}\omega^{-1}_j, \label{1.23}\\ 
 \omega_{N+1}+\omega_{N+2} &= -a_{N+1}-\sum_{j=1}^{N} \omega_j, \label{1.24}
 \end{align}
 where $\displaystyle a_0 = (-1)^{N}\omega^2_{\rm s}\prod_{j=1}^{N} t^{-1}_{j}$,
 $a_{N+1} = -\displaystyle\sum_{j=1}^{N}t^{-1}_{j} (1+\varepsilon_j^0).$
  
Equations \eqref{1.23} and \eqref{1.24} are equivalent to the second-order polynomial equation, which has two complex conjugate roots. By virtue of (\ref{eq:perturbation}), waves, which correspond to these roots, are dynamical and can damp. We take a wave, which corresponds to the value $\omega_{N+1}$ with a negative imaginary component as a reference solution. Through the perturbation of gas density $\delta\hat{f}$, we define all other perturbations by using the system \eqref{1.8} - \eqref{1.11} 

\begin{align}\label{1.13}
\frac{\delta\hat{v}}{c_s} &=-i\frac{\omega_{N+1}}{\omega_{\rm s}}\frac{\delta\hat{\rho}_{\rm g}}{\rho_{\rm g}^0},\\ \label{1.14}
\frac{\delta\hat{u}_j}{c_s}&=-i\frac{\omega_{N+1}}{\omega_{\rm s}}\frac{1}{(1-\omega_{N+1} t_{j})}\frac{\delta\hat{\rho}_{\rm g}}{\rho_{\rm g}^0},\\\ \label{1.15}
\frac{\delta\hat{\rho}_j}{\rho^0_j}&=\frac{1}{(1-\omega_{N+1} t_{j})}\frac{\delta\hat{\rho}_{\rm g}}{\rho_{\rm g}^0},
\end{align}
where $j=1,\cdots,N$.

Assume that $\delta\hat{\rho}_{\rm g}=A$, then the following solution in the plane of complex numbers will correspond to the wave with $\omega_{N+1}$: 
\begin{align}\label{1.16}
\delta\rho_{\rm g}(x,t) &=Ae^{ikx-\omega_{N+1} t},\\ \label{1.17}
\delta v(x,t) &=-Ai\frac{\omega_{N+1}}{\omega_{\rm s}}\frac{c_s}{\rho_{\rm g}^0}e^{ikx-\omega_{N+1} t},\\ \label{1.18}
\delta\rho_j (x,t)&=A\frac{\rho^0_j}{\rho_{\rm g}^0}\frac{1}{(1-\omega_{N+1} t_{j})}e^{ikx-\omega_{N+1} t},\\ \label{1.19}
\delta u_j (x,t)&=-Ai\frac{\omega_{N+1}}{\omega_{\rm s}}\frac{c_{\rm s}}{\rho_{\rm g}^0}\frac{1}{(1-\omega_{N+1} t_{j})}e^{ikx-\omega_{N+1} t}.
\end{align}

The real part of functions \eqref{1.16}-\eqref{1.19} satisfies the system (\ref{1.4})-(\ref{1.7}). We will use it as a reference solution. Then for $t=0$, the reference solution takes the form
\begin{equation}\label{1.27}
\delta f(x)=A[Re(\delta\hat{f})cos(kx)-Im(\delta\hat{f})sin(kx)].
\end{equation}

\section{Sound speed in two-phases polydisperse medium} 
\label{sec:analysis}
Phase velocity of the wave is defined as $v_{\rm ph}=\omega/k$. For sound waves in gas $\omega=c_{\rm s}k$ and  $v_{\rm ph}=c_{\rm s}$. In our case of the medium with polydisperse dust, the frequency is a complex function of wave vector. To find the sound speed, we need to analyze the expression  \eqref{1.12} for $\omega$. Relation \eqref{1.12} is an implicit dependence of frequency on the wave number and parameters of polydisperse medium with $N$ dust components:
\begin{equation}\label{1.29}
\omega^2(1+\sum_{j=1}^{N}\frac{\varepsilon^0_j}{1-\omega t_{j}})+k^2c_{\rm s}^2=0.
\end{equation}
We remind that relation \eqref{1.29} is obtained for the solutions of linearized system with perturbations of the form $\delta f=\delta\hat{f}e^{ikx-\omega t}$, i.e., the solution for sound wave corresponds to imaginary values $\omega$.\\
For low-frequency waves $\omega t_j\ll1$, the expression \eqref{1.29} is simplified to the form
\begin{equation}\label{1.30}
\omega^2(1+\sum_{j=1}^{N}\varepsilon^0_j)+k^2c_{\rm s}^2=0,
\end{equation}
From \eqref{1.30}, it follows that $\omega=\pm i kc_{\rm s}/\sqrt{1+\sum_{j=1}^{N}\varepsilon^0_j}$. Hence, we obtain that the sound speed in polydisperse medium within low-frequency limit 
\begin{equation}
\label{eq:csmulti}
    c^*_{\rm s}=v_{\rm ph}=c_{\rm s}/\sqrt{1+\sum_{j=1}^{N}\rho_j^0/\rho_{\rm g}^0}.
\end{equation}

\end{document}